\documentclass[a4paper,11pt,oneside,reqno]{amsart}
\usepackage{amscd,amsfonts,amssymb,amsmath,amsthm,bbm}
\usepackage{graphicx,psfrag}
\usepackage[english]{babel}
\usepackage[headinclude,DIV13]{typearea}
\areaset{15.1cm}{25.0cm}
\usepackage[hypertex,
  colorlinks=true,
 linkcolor=blue,
 citecolor=blue,
]{hyperref}
\usepackage[T1]{fontenc}
\usepackage{mathrsfs}
\usepackage{amsmath,amsthm,amssymb,psfrag}
\usepackage{latexsym}

\usepackage{enumerate}
\usepackage{mathrsfs}
\usepackage{amsopn}
\usepackage{amsmath}
\usepackage{amssymb}
\usepackage{amsfonts}
\usepackage{amsbsy}
\usepackage{amscd,indentfirst,epsfig}
\usepackage{hyperref}
\usepackage{amsfonts,amsmath,latexsym,amssymb,verbatim,amsbsy}
\usepackage{amsthm}
\renewcommand{\theequation}{\arabic{section}.\arabic{equation}}

\theoremstyle{plain}
\newtheorem{theo}{Theorem}[section]
\newtheorem{lemme}[theo]{Lemma}
\newtheorem{lemmeA}{Lemma A.}
\newtheorem{nbA}{Remark A.}

\newtheorem{propoA}{Proposition A.}

\newtheorem{propo}[theo]{Proposition}
\newtheorem{cor}[theo]{Corollary}

\newtheorem{nb}[theo]{Remark}
\newtheorem{defi}[theo]{Definition}
\theoremstyle{definition}

\def \leq {\leqslant}
\def \geq {\geqslant}

\numberwithin{equation}{section}
\def\ind#1{\lower5pt\hbox{$\scriptstyle #1$}}

\def \d {\,\mathrm{d} }
\def \L {\mathcal{L}}

\def \H {\mathcal{H}}

\def \ds {\displaystyle}

\def \LL {\mathscr{L}_2}

\def \ua {\mathbf{u}_\alpha}

\def \D {\mathscr{D}}

\def \ds {\displaystyle}
\def\Q {\mathcal{Q}}

\def\R{{\mathbb R}}

\def \S {{\mathbb S}^2}

\def \u {{u} }
\def \v {{v}}

\def \vb {\v_{\star}}
\def \w {{w}}

\def \fe {F_\alpha}
\def \he {G_\alpha}
\def \M {\mathcal{M}}
\def \sM {\mathsf{M}}

\def \IS {\int_{\S}}
\def \IR {\int_{\R^3}}
\def \IRR {\int_{\R^3 \times \R^3}}

\def \B {\mathcal{B}}
\def \G {\mathcal{G}}

\title[Unique steady state for hard-spheres driven by a particle bath]
{{Uniqueness in the weakly inelastic regime of the equilibrium state of the inelastic Boltzmann equation
driven by a particle  bath}}

\author{Marzia  Bisi, Jos\'{e}  A. Ca\~{n}izo \& Bertrand Lods }

\address{\textbf{Marzia Bisi}, Dipartimento di Matematica, Universit\`{a} di Parma,
    Viale G. P. Usberti 53/A, 43100 Parma, Italia}    \email{marzia.bisi@unipr.it}

\address{\textbf{Jos\'{e} A.  Ca\~{n}izo},  Departament de Ma\-te\-m\`a\-ti\-ques, Universitat
Aut\`onoma de Barcelona, E-08193 Bellaterra, Spain}
\email{canizo@mat.uab.cat}

\address{\textbf{Bertrand Lods}, Dipartimento di Statistica e Matematica Applicata \& Collegio Carlo Alberto, Universit\`{a} degli
Studi di Torino,  Corso Unione Sovietica, 218/bis, 10134 Torino, Italy.}\email{lods@econ.unito.it}

\begin{document}

\maketitle

\begin{abstract} We consider the spatially homogeneous Boltzmann equation for inelastic
hard-spheres (with constant restitution coefficient $\alpha \in (0,1)$)
 under the thermalization induced by a host medium with a fixed Maxwellian
  distribution. We prove uniqueness of the stationary solution (with given
  mass) in the weakly inelastic
  regime; i.e., for any inelasticity parameter $\alpha \in (\alpha_0,1)$,
  with some constructive $\alpha_0 \in  [0, 1)$.  Our analysis is based on
  a perturbative argument which uses the knowledge of the stationary solution in the elastic limit and quantitative
  estimates of the convergence of stationary solutions as the inelasticity parameter goes to $1$. In order to
  achieve this we give an accurate spectral analysis of the associated
  linearized collision operator in the elastic limit. Several qualitative
  properties of this unique steady state $\fe$ are also derived; in
  particular, we prove that $\fe$ is bounded from above and from below by two
  explicit universal (i.e. independent of $\alpha$) Maxwellian distributions.
\end{abstract}

\section{Introduction}
\setcounter{equation}{0}

\subsection{Physical context: driven granular gases} Kinetic models
for dilute granular flows are based, as well documented \cite{Br}, on
a Boltzmann equation in which collisions between hard--spheres
particles are supposed to be inelastic, i.e.\ at each encounter a
fraction of the kinetic energy is dissipated. Such a dissipation
implies that, in absence of energy supply, inelastic hard spheres are
cooling down and the energy continuously decreases in time. In
particular, the corresponding dissipative Boltzmann equation admits
only trivial equilibria.  This is no longer the case if the spheres
are forced to interact with an external agent (thermostat), in which
case the energy supply may lead to a non-trivial steady state. For
such driven system (in a space-homogeneous setting), the
time evolution of the one-particle distribution function \(f(v,t)\),
\(v\in\R^3\), \(t>0\) satisfies the following
\begin{equation}
\label{be:force} \partial_t f  =\Q_\alpha(f,f) + \G(f),
\end{equation}
where \(\Q_\alpha(f,f)\) is the inelastic quadratic Boltzmann
collision operator (see next section for a precise definition), while
\(\G(f)\) models the forcing term. The parameter $\alpha \in (0,1)$ is
the so-called ``restitution coefficient'', expressing the degree of
inelasticity of binary collisions among grains and the purely elastic
case is recovered when $\alpha = 1$.

There exist in the literature several possible physically meaningful
choices for the forcing term $\G$ in order to avoid the cooling of the
granular gas. The most natural one is the pure diffusion thermal bath
for which the particles are subject to uncorrelated random
accelerations between the collisions yielding to the diffusive
operator
$$\G_{1} (f) = \mu \,\Delta {f},$$
where $\mu >0$ is a given parameter and $\Delta=\Delta_v$ the
Laplacian in the velocity variable. For this model, introduced in
\cite{noije}, the existence of a non-trivial equilibrium state has
been obtained in \cite{GPV**} while the uniqueness (in some weakly
inelastic regime) and the linear/nonlinear stability of such a steady
state has been proved in \cite{MiMo3}.  Other fundamental examples of
forcing terms are the thermal bath with linear friction \cite{Barrat}:
$ \G_2(f)= \lambda \,\Delta {f} + \kappa\, \mathrm{div} (v\,f)$ with
several range of parameters $\kappa,\lambda$ and where $\mathrm{div}$ is the
divergence operator with respect to the velocity variable. A
particular case of interest is the one related to the following
anti-drift operator
$$\G_3 (f) = -\mu \,\mathrm{div} (v f), \qquad \mu >0.$$
The existence of an equilibrium state for such a forcing term has been proved in \cite{MiMo2} and is related to the existence of self-similar solutions to the freely evolving Boltzmann equation. Such a steady state corresponds then to a self-similar profile (the so-called Homogeneous Cooling State) and both its uniqueness  (still in some weakly inelastic regime) and its stability have been derived in \cite{MiMo2}, providing a rigorous proof to  the Ernst-Brito
conjecture~\cite{ernstbrito} for inelastic hard-spheres in the weakly inelastic regime.  \medskip

\subsection{Description of the problem and main results} In this paper
we address a problem similar to the aforementioned ones but with a
forcing term of different nature.  Namely, we consider a situation in
which the system of inelastic hard spheres is immersed into \textit{a
  thermal bath of particles} so that the forcing term $\G$ is given by
a linear scattering operator describing inelastic collisions with the
background medium.  More explicitly, we shall assume in the present
paper that the forcing operator $\G$ is a linear Boltzmann collision
operator of the form:
$$\G(f)=: \L (f) = \Q_e (f, \M_0)$$
where $\Q_e(\cdot,\cdot)$ is a Boltzmann collision operator associated
to the (fixed) restitution coefficient $e \in (0,1]$ and $\M_0$ stands
for the distribution function of the host fluid. We shall assume here
that this host distribution is a given Maxwellian with unit mass, bulk
velocity~$u_0$ and temperature $\Theta_0>0$:
\begin{equation}\label{maxwe1}
\M_0(\v)=\bigg(\dfrac{1}{2\pi
\Theta_0 }\bigg)^{3/2}\exp
\left\{-\dfrac{(\v-\u_0 )^2}{2\Theta_0 }\right\}, \qquad \qquad \v
\in \R^3.
\end{equation}
The precise definitions of both collision operators $\Q_\alpha(f,f)$ and $\L (f)$,
with their weak forms and the relations between pre- and post--collision velocities,
are given in \textsc{Subsection~2.1}.\smallskip

The existence of smooth stationary solutions for the inelastic
Boltzmann equation under the thermalization induced by a host-medium
with a fixed distribution has been investigated by two of the authors,
in collaboration with J. A. Carrillo, in \cite{BiCaLo}; we refer to
this paper the references therein for more information about the
physical relevance of such thermal bath of particles. To be more
precise, it has been proved in \cite{BiCaLo} that, for any restitution
coefficient $\alpha \in (0,1]$, there exists a non-trivial smooth
stationary state $F_\alpha \geq 0$ such that
\begin{equation}\label{stead}\Q_\alpha(F_\alpha,F_\alpha)+\L(F_\alpha)=0.\end{equation}
The proof of this existence result is based on a dynamic version of
Tykhonov fixed point theorem and is achieved by controlling the
$L^p$-norms, the moments and the regularity of the solutions for the
Cauchy problem \eqref{be:force}. Moreover, using the analysis of the
linear scattering operator $\L$, for elastic nonlinear interactions
(i.e. whenever $\alpha=1$) one can prove easily that there exists a \textit{unique solution} with unit mass to the equation
  \begin{equation}\label{equi1}
    \Q_1(F,F)+\L(F)=0.\end{equation}
Moreover, this unique distribution is a Maxwellian $\M(v)$ with bulk velocity $u_0$ and explicit temperature $\Theta^{\#} \leq \Theta_0.$ The knowledge of the equilibrium solution in the elastic case
$\alpha=1$ will be of paramount importance in our analysis of the
steady state $\fe$ in the weakly inelastic regime $\alpha \simeq 1.$
All these preliminary results are recalled in
\textsc{Section~2}. \medskip

Uniqueness and qualitative properties of the steady distribution are
still open problems for $\alpha < 1$, and these are the main subjects
of the present paper. To be more precise, as far as uniqueness is
concerned, we prove the following:
\begin{theo}\label{mainUni}
  There exists $\alpha_0 \in (0,1]$ such that, for any $\varrho \geq 0$, the set
  $$\mathfrak{S}_\alpha(\varrho)=\left\{F_\alpha \in L^1_2, \:F_\alpha \geq 0,\; \fe  \text{ solution to \eqref{stead} with } \IR \fe(v) \d v =\varrho  \right\}$$
  reduces to a singleton where $L_2^1$ is the set of integrable
  distributions with finite energy. In particular, for any $\alpha \in
  (\alpha_0,1]$, such a steady state $F_\alpha$ is radially symmetric
  and belongs to $\mathcal{C}^\infty(\R^3)$.
\end{theo}

Several further qualitative properties of the steady state $\fe$ are
also given in the paper.  In particular, we are able to derive
pointwise estimates for the steady state $\fe$ which are uniform with
respect to the inelasticity parameter $\alpha$:

\begin{theo}\label{pointIn}
  There exist two Maxwellian distributions $\underline{\M}$ and $\overline{\M}$
  (independent of $\alpha$) such that
  \begin{equation}
    \label{pointwise}
    \underline{\M}(v) \leq \fe(v) \leq \overline{\M}(v)
    \qquad \forall v \in \R^3,
    \qquad \forall \alpha \in (0,1).
  \end{equation}
\end{theo}
For the upper bound, the strategy of proof is inspired by the
comparison principle of~\cite{GaPaViM} and uses some estimates of
\cite{AlCaGa}. For the lower Maxwellian bound, the proof is much
simpler than the ones yielding (non Maxwellian) pointwise lower bounds
for the forcing terms $\G_1$ and $\G_3$ (see \cite{MiMo2,MiMo3}) which
rely on the spreading properties of the quadratic inelastic collision
operator $\Q_\alpha$. Our approach relies uniquely on the properties
of the linear collision operator $\L$ and, more precisely, on the
explicit integral representation of the gain part $\L^+$ derived in
\cite{ArLo}.\medskip

More general theorems analyzing possible global stability properties
of the stationary solution are planned as future work.

\subsection{Strategy of proof and organization of the paper}

Our strategy of proof is inspired by the strategies adopted in
Refs.~\cite{MiMo2, MiMo3} for different kinds of forcing
terms. However, the peculiarities of our linear scattering operator
such as its lack of symmetry and the exchange of momentum between
grains and background, will require in some points a completely
different treatment with respect to previous works on analogous
problems.

\subsubsection{Main difference with respect to other forcing terms}

Let us spend a few words in explaining the key differences (which will
also be emphasized throughout the paper):
\begin{itemize}
\item The quadratic operator $\Q_\alpha(f,f)$ preserves mass and
  momentum and both the forcing terms $\G_1$ and $\G_3$ considered in
  Refs.~\cite{MiMo2, MiMo3} also do so.  Therefore, for both these
  forcing terms, the mass and momentum of a stationary solution can be
  prescribed. This is no more the case whenever the forcing term is
  the linear scattering $\L$ which \textit{does not preserve
    momentum}.
\item Moreover, while the collisional operator $\Q_\alpha$ tends to
  cool down the gas --- dissipating kinetic energy --- the forcing
  terms $\G_1$ and $\G_3$ have the tendency to warm it up in some explicit way. Precisely, for any nonnegative
  distribution $f$,
  $$\IR \G_1(f)\,|v|^2\,\d v=6\mu\varrho_f
  \qquad \text{ while } \qquad \IR \G_3(f)\,|v|^2\d v=2\mu
  \mathcal{E}_f$$
  where $\varrho_f=\IR f(v)\d v$ is the prescribed mass density of
  $f$ and $\mathcal{E}_f=\IR |v|^2 f(v)\d v$ denotes its energy. It is unfortunately
  impossible to quantify the thermal contribution of the linear
  scattering operator $\L$ in such a closed way: indeed, since we are
  dealing with a linear scattering operator associated to hard-spheres
  interactions, the thermal contribution $\IR \L(f)\,|v|^2\d v$
  involves moments of $f$ up to third order.
\item Finally, it is not possible in our case to use the fundamental
  scaling argument of~\cite{MiMo2,MiMo3}. Precisely, for the forcing
  terms $\G_1$ and $\G_3$ studied in~\cite{MiMo2,MiMo3}, scaling
  arguments show that it is possible to choose $\mu >0$ arbitrarily and this
  yields the authors of~\cite{MiMo2,MiMo3} to choose $\mu=\mu_\alpha$
  so that, in the elastic limit $\alpha \to 1$, the dissipation of
  kinetic energy will exactly be balanced by the forcing term. Such a
  scaling argument \textit{cannot be invoked} for the linear
  scattering operator $\L$ and this is again related to the fact that
  we are dealing here with hard-spheres interactions. Notice that, if
  $\L$ were the linear Boltzmann operator associated to
  pseudo-Maxwellian molecules, the scaling argument and, more
  generally, the whole strategy of \cite{MiMo2,MiMo3} would apply
  almost directly (still assuming $\Q_\alpha$ to be associated to
  hard-spheres).
\end{itemize}

\subsubsection{General strategy}

Let us now explain the main steps in our strategy of proof. It is
essentially based on the knowledge of the elastic limit problem and on
quantitative estimates of the difference between solutions to the
original problem and the equilibrium state in the elastic limit. Introduce the linearized operator in the elastic limit (where we recall that, for $\alpha=1$, the unique steady state is an explicit Maxwellian $\M$)
\begin{equation}
  \label{linearized1}
  \mathscr{L}_1 (h)=\Q_1(\M,h)+\Q_1(h,\M) +\L h
\end{equation}
Given $\alpha \in (0,1]$, let $\fe$ and $\he$ belong to $\mathfrak{S}_\alpha(\varrho)$. One has
$$\Q_\alpha(\fe,\fe)+\L(\fe)=0=\Q_\alpha(\he,\he)+\L(\he).$$
Then
$$\mathscr{L}_1(\fe-\he)=\Q_1(\M,\fe-\he)+\Q_1(\fe-\he,\M)-\Q_\alpha(\fe,\he)+\Q_\alpha(\he,\he).$$
It is easy to recognize then that
\begin{multline*}
\mathscr{L}_1(\fe-\he)=\bigg(\Q_1(\fe-\he,\M) -\Q_\alpha(\fe-\he,\M)\bigg) \\
+\bigg(\Q_1(\M,\fe-\he) -\Q_\alpha(\M,\fe-\he)\bigg)\\
+\bigg(\Q_\alpha(\fe-\he,\M-\fe) -\Q_\alpha(\M-\he,\fe-\he)\bigg).
\end{multline*}
Assume now that there exist two Banach spaces $\mathcal{X}$ and
$\mathcal{Y}$ (independent of $\alpha$) such that
\begin{multline}
  \label{XY}
  \left\|\Q_1(h,\M)-\Q_\alpha(h,\M)\right\|_{\mathcal{X}}
  \\
  +\left\|\Q_1(\M,h)-\Q_\alpha(\M,h)\right\|_{\mathcal{X}}
  \leq
  \eta(\alpha) \|h\|_\mathcal{Y},
  \qquad \forall \alpha \in (0,1)
\end{multline}
where $\lim_{\alpha \to 1}\eta(\alpha)=0$ and there exists
$C >0$ such that
\begin{equation}
  \label{XY1}
  \left\|\Q_\alpha(h,g)\right\|_{\mathcal{X}}
  +\left\|\Q_\alpha(g,h)\right\|_{\mathcal{X}}
  \leq
  C \,\|g\|_\mathcal{Y} \|h\|_\mathcal{Y},
  \qquad \forall \alpha \in (0,1),
\end{equation}
then,
\begin{multline*}\left\|\mathscr{L}_1(\fe-\he)\right\|_\mathcal{X} \leq \left(\eta(\alpha)+  C \left\|\fe-\M\right\|_\mathcal{Y} +\right. \\
+\left. C \left\|\he-\M\right\|_\mathcal{Y}\right)\left\|\fe - \he\right\|_\mathcal{Y} \qquad \forall \alpha \in (0,1).\end{multline*}
If moreover there exists $c_0 >0$ such that
\begin{equation}\label{XY2}
\|\mathscr{L}_1 (h)\|_\mathcal{X} \geq c_0 \|h\|_\mathcal{Y} \qquad \forall h \in \bigcup_{\alpha \in (0,1)}\mathfrak{S}_\alpha(0) \subset \mathcal{Y}\,,
\end{equation}
then one sees  that
$$c_0\left\|  \fe-\he \right\|_\mathcal{Y} \leq \delta(\alpha)
\left\|\fe - \he\right\|_\mathcal{Y} \qquad \forall \alpha \in (0,1)$$
with $$\delta(\alpha)=\eta(\alpha)+  2C \max\left\{\left\|\fe-\M\right\|_\mathcal{Y},\left\|\he-\M\right\|_\mathcal{Y}\right\}.$$
Therefore, if we are able to construct an explicit $\alpha_0 \in (0,1)$ such that
\begin{equation}\label{deltac0}
\fe, \he \in \mathfrak{S}_\alpha \text{ with } \alpha \in (\alpha_0,1] \implies \delta(\alpha) < 1/c_0
\end{equation}
then
$$\fe=\he.$$
All the technical difficulty is then to determine $\mathcal{X}$ and $\mathcal{Y}$ such that \eqref{XY}, \eqref{XY1} and \eqref{XY2} hold true and to prove that the $\mathcal{Y}$-norm is compatible with \eqref{deltac0}:
\begin{enumerate}
\item The proof that \eqref{XY} and \eqref{XY1} hold true will be straightforward on the basis of known estimates of the collision operator $\Q_\alpha$.
\item Notice that \eqref{XY2} means that $\mathscr{L}_1 \::\:\bigcup_{\alpha \in (0,1)}\mathfrak{S}_\alpha(0) \subset \mathcal{Y} \to \mathcal{X}$ is invertible and the proof of such a property relies on a careful spectral analysis of $\mathscr{L}_1$.
\item Concerning now estimate \eqref{deltac0}, it consists in proving that
  $$\lim_{\alpha \to 1} \sup_{\fe \in \mathfrak{S}_\alpha}\left\|\fe-\M\right\|_{\mathcal{Y}}=0$$
  More precisely, it amounts to providing a quantitative estimate on
  the distance between $\fe$ and the Maxwellian $\M$ in the elastic
  limit $\alpha \to 1$. This is the most technical part of the
  uniqueness result.
\end{enumerate}

To be able to complete the above program, one begins with deriving
suitable \emph{a posteriori} estimates on the steady state that shall
be useful in the sequel.  In particular, after estimating the
high--energy tails of the solution $f(v,t)$ to the Boltzmann equation
\eqref{be:force} uniformly with respect to the inelasticity parameter
$\alpha$, it is possible to prove that, for any $\alpha \in (0,1]$,
the stationary solution $F_\alpha$ admits an exponential tail of
second-order. Moreover, we obtain uniform lower and upper bounds on
the energy of $F_\alpha$ and this yields a control of $H^k$--norms.

To prove the points (2) and (3) of the above program, we derive the
spectral properties of the linearized collision operator in the
elastic limit $\mathscr{L}_1$ given by \eqref{linearized1}. As already
mentioned, this quantitative spectral analysis of $\mathscr{L}_1$
resorts to very recent results \cite{gualdani} which allow to extend a
spectral gap result from a smaller (typically Hilbert) space
$\mathcal{H}$ to a larger (typically Banach) space $\mathcal{X}$. We
apply these recent abstract results to both the linear scattering
operator $\L$ (whose spectral analysis in a weighted $L^2$-space has
been performed in \cite{ArLo,LoMo}) and to the linearized operator
$\mathscr{L}_1$. This will allow us to prove point (2) of the above
program. Moreover, this spectral analysis will also allow to provide a
quantitative estimate on the distance between~$\fe$ and the
Maxwellian~$\M$ in the quasi elastic limit $\alpha \to 1$ (see point
(3) above). We wish to emphasize here the fact that, with our
approach, we prove the convergence of $\fe$ to $\M$ as $\alpha \to 1$
without knowing \textit{a priori} that the energy $\mathcal{E}_{\fe}$
converges to that $\mathcal{E}_{\M}$ of the Maxwellian $\M$. This is a
major difference with respect to the papers \cite{MiMo2,MiMo3} where,
for the reasons already explained, it was possible to write down a
relatively simple equation (in closed form) satisfied by the
difference $\mathcal{E}_{\fe}-\mathcal{E}_{\M}$. This is not possible
in the present situation since, again, $\L$ is a scattering operator
associated to hard-spheres interactions.

\subsubsection{Organization of the paper} After recalling the precise definitions of the Boltzmann operator $\Q_\alpha$
and the forcing term $\L$, we give a precise simple proof of uniqueness of the equilibrium in the elastic case in \textsc{Section 2}.
Then,  \textsc{Section 3} is devoted to the derivation of the  \emph{a posteriori} estimates on the steady state for general restitution coefficients.
The uniform pointwise estimates (Theorem \ref{pointIn}) are proved in \textsc{Section~4} while \textsc{Section 5} is devoted to the proof
of Theorem~\ref{mainUni}. In an Appendix of the paper, several estimates on $\L$ and $\mathscr{L}_1$ are derived which turn out to be useful for the spectral analyis performed in \textsc{Section 5}.

\section{Preliminary results}
\subsection{The kinetic model}
Given a constant restitution coefficient $\alpha \in (0,1)$, one
defines the bilinear Boltzmann operator $\Q_\alpha$ for inelastic
interactions and hard-spheres by its action on test functions
$\psi(v)$:
 \begin{equation}\label{weakfg}
 \IR \Q_\alpha (f,g)(v)\, \psi(v)\d\v  = \dfrac{1}{ 4\pi}\IR\IR \IS f(v)g(w)\,|v-w|\,\left(\psi(v')-\psi(v)\right)\d v \d w \d \sigma
 \end{equation}
 with $v'=v+\frac{1+\alpha}{4}\,(|v-w|\sigma-v+w)$. In particular, for
 any test function $\psi=\psi(v)$, one has the following weak form of
 the \textit{quadratic} collision operator:
 \begin{equation}\label{co:weak}
\begin{split}
\IR \Q_\alpha (f,f)(v)\, \psi(v)\d\v  =  \frac{1}{2} \IR\IR f(v)\,f(w)\,|v-w|
\mathcal{A}_\alpha[\psi](v,w)\d\w\d\v,
\end{split}
\end{equation}
where
\begin{equation}\begin{split}
\label{coll:psi} \mathcal{A}_\alpha[\psi](v,w) &=
\frac{1}{4\pi}\IS(\psi(v')+\psi(w')-\psi(v)-\psi(w))\d{\sigma}\\
&=\mathcal{A}^+_\alpha[\psi](v,w)-\mathcal{A}^-_\alpha[\psi](v,w) \end{split}
\end{equation}
and the post-collisional velocities $(v',w')$ are given by
\begin{equation}
\label{co:transf}
  v'=v+\frac{1+\alpha}{4}\,(|q|\sigma-q),
\qquad  w'=w-\frac{1+\alpha}{4}\,(|q|\sigma-q), \qquad q=v-w.
\end{equation}
In the same way, for another constant restitution coefficient $e \in
(0,1)$, one defines the linear scattering operator $\L$ by its action
on test functions:
 \begin{equation}\label{co:weak}
\begin{split}
\IR \L (f)(v)\, \psi(v)\d\v  = \IR\IR f(v)\,\M_0(w)\,|v-w|
\mathcal{J}_e[\psi](v,w)\d\w\d\v,
\end{split}
\end{equation}
where
\begin{equation}\label{col:Je}
 \mathcal{J}_e[\psi](v,w)=
\dfrac{1}{4\pi}\int_{\S}
\left(\psi( {v}^\star)-\psi(v)\right)\d\sigma=\mathcal{J}_e^+[\psi](v,w)-\mathcal{J}_e^-[\psi](v,w).
\end{equation}
with  post-collisional velocities $(v^\star,w^\star)$
\begin{equation}
\label{co:transf1}
  v^\star=v+\frac{1+e}{4}\,(|q|\sigma-q),
\qquad  w^\star=w-\frac{1+e}{4}\,(|q|\sigma-q), \qquad q=v-w.
\end{equation}
For simplicity, we shall assume in the paper that the total mass of
the particles governed by $f$ and that of $\M_0$ are equal. Notice
that
$$\L(f)=\Q_e(f,\M_0)$$
and we shall adopt the convention that post (or pre-) collisional
velocities associated to the coefficient $\alpha$ are denoted
with \textit{prime} symbol, while that associated to $e$ are denoted
with $\star$ symbol. We are interested in the stationary solution to
the following Boltzmann equation:
\begin{equation}
  \label{BEev}
  \partial_t f(t,v)
  = \Q_\alpha(f(t,\cdot);f(t,\cdot))(v) + \L(f)(t,v)
  \qquad t > 0, \quad f(0,v)=f_0(v).
\end{equation}
We proved the following in \cite{BiCaLo}
\begin{theo}[\textit{\textbf{Existence of stationary solutions}}]
  \label{theo:main}
  For any restitution coefficient $\alpha \in (0,1)$, there exists a nonnegative $F_\alpha \in
  L^1_2\cap L^p$, $p\in (1,\infty)$ with unit mass and positive
  temperature such that
\begin{equation}\label{solutionstationaire}\Q_\alpha(F_\alpha,F_\alpha) +\L(F_\alpha)=0.\end{equation}
  Moreover, there exists a steady state which is radially symmetric and belongs to $\mathcal{C}^\infty(\R^3)$.
\end{theo}
\begin{nb} Notice that the existence of a radially symmetric
  stationary solution to \eqref{solutionstationaire} is not explicitly
  stated in \cite{BiCaLo}, where more general host distributions than
  $\M_0$ are considered. However, since the Maxwellian distribution
  $\M_0$ is radially symmetric, one easily checks that the property of
  being radially symmetric is stable along the flow of \eqref{BEev}
  (i.e. $f_0$ radially symmetric $\implies f(t,v)$ radially symmetric
  for any $t\geq 0$) and therefore the fixed point argument used in
  \cite{BiCaLo} allows to build a radially symmetric steady solution
  to \eqref{solutionstationaire}.
\end{nb}

Notice that
$$\Q_\alpha(f,f)=\Q^+_\alpha(f,f)-\Q^-_\alpha(f,f)=\Q^+_\alpha(f,f)-f \Sigma(f) $$
where
$$\Sigma(f)(v)=(f \ast |\cdot|)(v)=\IR f(w)|v-w|\d w.$$
Notice that $\Sigma(f)$ does not depend on the restitution coefficient $\alpha \in (0,1]$. In the same way,
$$\L(f)(v)=\L^+(f)(v)-\L^-(f)(v)=\L^+(f)(v)-\sigma(v)f(v)$$
where
\begin{equation}\label{sigma}
\sigma(v)=(\M_0 \ast |\cdot|)(v)=\IR \M_0(w)|v-w|\d w.\end{equation}

\subsection{A basic observation in the elastic case} We begin with a
basic observation concerning the elastic case. Precisely, when the
quadratic operator is that for elastic interactions, i.e. for
$\alpha=1$, one can prove in a very direct way that the steady state
solution to the above problem is unique. Precisely, the background
forces the system to adopt a Maxwellian steady state (with density
equal to $1$):
\begin{theo}\label{steadystate} The Maxwellian velocity distribution:
\begin{equation}\mathcal{M} (\v)=\left(\dfrac{1}{2\pi
\Theta^{\#}}\right)^{3/2}\exp
\left\{-\dfrac{(\v-\u_0)^2}{2\Theta^{\#}}\right\}, \qquad \v \in
\R^3,\end{equation} with
\begin{equation}
  \Theta^{\#}
  =\dfrac{1+e}{3-e}\Theta_0
\end{equation}
is the unique solution with unit mass to the equation
\begin{equation}\label{equi}
  \Q_1(F,F)+\L(F)=0.\end{equation}
\end{theo}

\begin{proof}
  It has been proved in \cite{LoTo} that $\L(\mathcal{M})=0$. Now,
  since $\mathcal{M}$ is a Maxwellian distribution, it is also
  well-known that $\Q_1(\mathcal{M},\mathcal{M})=0$ and this proves
  that $\mathcal{M}$ is a solution to~\eqref{equi}. To prove that it
  is the unique solution with unit mass, one proceeds {in some formal way for the time being assuming that $F$ decays sufficiently fast at infinity; we will see that it is actually the case in the following section. All the proof can then be made rigorous thanks to the subsequent Theorem \ref{gaus}}. For
  any distribution $F(\v) \geq 0$ solution to \eqref{equi}, let us
  multiply \eqref{equi} with $\log
  \left(\frac{F(\v)}{\mathcal{M}(\v)}\right)$ and integrate with
  respect to $\v$. One gets
  $$0=\IR \Q_1(F,F)(\v)\log
  \left(\frac{F(\v)}{\mathcal{M}(\v)}\right)\d\v + \IR \L(F)(\v)\log
  \left(\frac{F(\v)}{\mathcal{M}(\v)}\right)\d\v$$ and it is
  well-known from \cite{LoTo} and \cite{Ce94} that both the integrals
  in the above sum are nonpositive. Therefore,
  $$\IR \Q_1(F,F)(\v)\log
  \left(\frac{F(\v)}{\mathcal{M}(\v)}\right)\d\v =0.$$ Since
  $\mathcal{M}$ is a Maxwellian distribution, it is a well-established
  fact that
  $$\IR \Q_1(F,F)(\v) \log \mathcal{M}(\v) \d \v=0.$$
  Consequently, $F$ is such that
  $$\IR \Q_1(F,F)(\v)\log F(\v) \d \v=0$$ and the classical Boltzmann
  $H$-Theorem \cite{Ce94} asserts that $F$ is a given Maxwellian and
  $\Q_1(F,F)=0.$ Consequently, one has $\L(F)=0$ and, from the
  uniqueness result \cite{LoTo}, $F=\mathcal{M}.$
\end{proof}

\section{A posteriori estimates}

\subsection{High-energy tails for the steady solution}\label{sect:high}

We are interested here in estimating the high-energy tails both of the solution $f(t,v)$  to \eqref{BEev} and of the stationary solutions to~\eqref{equi}  through a  weighted integral bound. Our approach is reminiscent to the work of \cite{Bobmom} recently improved in a series of papers \cite{BoGaPa, MiMo2, AloGam, GaPaViM, AloLo}.
\begin{defi}We say that the function $f$ has an {\it exponential
tail of order} $s>0$ if the following supremum
\begin{equation}
\label{tail_temp}
r^*_s = \sup \big\{r>0 \,|\,\mathcal{F}_{r,s}(f):=\int_{\R^3}f(v)\exp(r|v|^s)\d\v < +\infty \big\}
\end{equation}
is positive and finite.
\end{defi}
We begin by showing that, for the solution to \eqref{BEev}, exponential tails of order $s$ propagate with time if $s \in (0,2]$. The proof is adapted from several known results and follows the lines of \cite[Section 6]{AloLo}.
\begin{theo}\label{tailT}
Let $f_0$ be a nonnegative  velocity function with $\IR f_0(v)\d v=1$. Assume that $f_0$ has an exponential tail of order $s \in (0,2]$, i.e. there exists $r_0 >0$ and $s \in (0,2]$ such that
\begin{equation*}
\IR f_0(v)\exp\left(r_0|v|^s\right) \d v <\infty.
\end{equation*}
Then, there exist $0 < r \leq r_0$ and $C >0$ (independent of $\alpha \in (0,1]$) such that the  solution $f_\alpha(t,v)$ to the Boltzmann equation
\begin{equation}
  \label{BEevB}
  \partial_t f(t,v)
  = \Q_\alpha(f(t,\cdot);f(t,\cdot))(v) + \L(f)(t,v)
  \qquad t > 0, \quad f(0,v)=f_0(v)
\end{equation} satisfies
\begin{equation}\label{expbounds}
\sup_{t \geq0}\IR f_\alpha(t,v)\exp\left(r|v|^s\right) \d v \leq C < \infty.
\end{equation}
\end{theo}
\begin{proof} We adapt the strategy of \cite{BoGaPa} following carefully the dynamical approach of \cite{GaPaViM,AloLo}.
For notations convenience, we shall drop the dependence on $\alpha$ for the solution to \eqref{BEevB} and simply denote by $f(t,v)$ its solution.
Recall that, formally, $$\IR f(t,v)\exp\left(r|v|^s\right) \d v=\sum^{\infty}_{k=0}\frac{r^{k}}{k!}{m}_{\frac{sk}{2}}(t)$$
where
$${m}_p(t)=\IR f(t,v)|v|^{2p}\d v \qquad \forall t \geq 0, p \geq 1.$$
Therefore, to prove the result, it is sufficient to prove that there exists some $r >0$ (independent of $\alpha$) such that
$$\sum^{\infty}_{k=0}\frac{r^{k}}{k!}{m}_{sk/2}(t)
\text{ converges for any } t \geq 0.$$
From the Cauchy-Hadamard formula giving the radius of convergence of a
power series, it is enough to prove that, for any $s \in (0,2]$, there
is some real number $C=C(s) >0$ (independent of $t$ and of $\alpha$)
such that
\begin{equation}\label{Hadamard}
{m}_{\frac{sk}{2}}(t)  \leq  k! C^k  \qquad \forall t \geq 0, \:\forall k \in \mathbb{N}.\end{equation}
It is clear that, for any $p \geq 1$, the evolution of the $p$-moment ${m}_p(t)$ is given by
$$\dfrac{\d}{\d t}m_p(t)= {Q}_p(t)+ L_p(t)$$
where
$${Q}_p(t)=\IR \Q_\alpha(f(t,\cdot),f(t,\cdot))(v) |v|^{2p}\d v \quad \text{ and } \quad L_p(t)=\IR \L(f)(t,v)|v|^{2p}\d v.$$
Recall that the weak form of $\Q_\alpha$ and $\L$ are given in \eqref{co:weak}, \eqref{coll:psi} and \eqref{col:Je}. Now, based upon a sharp version of Povzner's estimates,  Bobylev, Gamba and Panferov
\cite[Lemma 1 \& Corollary 1]{BoGaPa} proved that, for any  $p \geq 1$,
$$\mathcal{A}^+_\alpha[|\cdot|^{2p}](v,w) \leq \gamma_{\alpha,p}\left(|v|^2+|w|^2\right)^p$$
where, for any $p >1$, $$\gamma_{\alpha,p}=\int_{-1}^1 \left(\frac{1+x}{2}\right)^p h_\alpha(x)\d x$$ with $h_\alpha(x)=\frac{1}{2}\left(g_\alpha(x)+g_\alpha(-x)\right)$ and
$$g_\alpha(x)=\dfrac{\left((1-\alpha)x+\sqrt{(1-\alpha)^2 x^2+4\alpha}\right)^2}{(1+\alpha)\sqrt{(1-\alpha)^2 x^2
+4\alpha}}, \qquad \forall x \in (-1,1).$$
Since we are looking for estimates which are uniform with respect to the inelasticity parameter $\alpha$, one notices that, as pointed out in \cite{MiMo,MiMo2},
for any $p\geq 1$, $$\sup_{\alpha \in (0,1)}\gamma_{\alpha,p} < \gamma_p:=\min\left(1,\frac{4}{p+1}\right).$$
In the same way,  one obtains the following very rough estimate for $\mathcal{J}^+_e[|\cdot|^{2p}](v,w)$:
$$\mathcal{J}_e^+[|\cdot|^{2p}](v,w) \leq \dfrac{1}{4\pi}\int_{\S}
\Big( |{v}^\star|^{2p}+|w^\star|^{2p} \Big) \d\sigma   \leq \gamma_{ p}\left(|v|^2+|w|^2\right)^p.$$ In particular, one obtains the following bounds:
\begin{multline}\label{Ae}
\mathcal{A}_\alpha[|\cdot|^{2p}](v,w) \leq \gamma_{ p}\left(|v|^2+|w|^2\right)^p -|v|^{2p}-|w|^{2p}\\
=-(1-\gamma_{ p})\left(|v|^{2p}+|w|^{2p}\right) + \gamma_{ p}\Big( (|v|^2+|w|^2)^p - |v|^{2p}-|w|^{2p} \Big)\end{multline}
and
\begin{multline}
  \label{Ie}
  \mathcal{J}_e[|\cdot|^{2p}](v,w)  \leq \gamma_{p}\left(|v|^2+|w|^2\right)^p -|v|^{2p}\\
  =-(1-\gamma_{ p}) |v|^{2p} + \gamma_{ p}\Big( (|v|^2+|w|^2)^p -
  |v|^{2p}-|w|^{2p}\Big)+\gamma_{ p}|w|^{2p}.\end{multline} Then, as
in \cite[Lemma 2 \& Eq. (4.5)]{BoGaPa}, one sees that
\begin{multline}
  \label{bo45}
  |v-w| \Big[ \left(|v|^2+|w|^2\right)^p -|v|^{2p}-|w|^{2p} \Big]
  \\
  \leq \sum_{k=1}^{k_p} \binom{p}{k}
  \left(|v|^{2(k+1/2)}|w|^{2(p-k)} +|v|^{2(p-k+1/2)}|w|^{2k}\right)
\end{multline}
where $k_p=\left[\frac{p+1}{2}\right]$ is the integer part of
$\frac{p+1}{2}$.  Performing now the $v$ and $w$ integrations and
using Eqs. \eqref{Ie}--\eqref{bo45} we get
\begin{equation*}
  \begin{split}
    L_p(t)  &\leq \frac{\gamma_{ p}}{2}\sum_{k=1}^{k_p}
    \binom{p}{k}
    \left(m_{ k+1/2 }(t)M_{ p-k }+m_{ p-k+1/2 }(t) M_{ k}\right)\\
    & \phantom{++++} +   \gamma_{ p}\left(m_{1/2}(t)M_p +m_0(t) M_{p+1/2}\right)\\
    &\phantom{+++++++} - \frac{1-\gamma_{ p}}{2}\IRR f(t,v)\M_0(w)|v-w|\,|v|^{2p}\d\v\d\w.
  \end{split}
\end{equation*}
One estimates the loss term using Jensen's inequality (together with
the fact that $u_0=0$) to get
$$ \IRR f(t,v)\M_0(w)|v-w|\,|v|^{2p}\d\v\d\w \geq \IR f(t,v)|v|^{2p+1}\d\v=m_{p+1/2}(t)$$
and
\begin{equation}
  \label{LPt}
  L_p(t) \leq -\, \frac{1-\gamma_{p}}{2} m_{p+1/2}(t)+     \frac{\gamma_{ p}}{2}\left(m_{1/2}(t)M_p +m_0(t) M_{p+1/2}\right)
  + \gamma_{ p}  \widetilde{S}_p(t)
\end{equation}
with
$$\widetilde{S}_p(t)=\frac{1}{2} \sum_{k=1}^{k_p}
\binom{p}{k} \left(m_{ k+1/2 }(t)M_{ p-k }+m_{ p-k+1/2 }(t) M_{ k}\right).$$
In the same way,  using Eqs. \eqref{Ae} and \eqref{bo45}, the following estimate was derived in \cite[Lemma 3]{BoGaPa}:
\begin{equation}\label{QPt}
Q_p(t) \leq -(1-\gamma_{ p})m_{p+1/2}(t)+ \gamma_{ p}   S_p(t)
\end{equation}
where
$$S_p(t)=\frac{1}{2}\sum_{k=1}^{k_p}
\binom{p}{k} \left(m_{ k+1/2 }(t)m_{ p-k }(t)+m_{ p-k+1/2 }(t) m_{ k}(t)\right).$$
This yields to the following differential inequality for $m_p(t):$
\begin{multline*}
\dfrac{\d}{\d t}m_p(t) \leq -\frac{3(1-\gamma_p)}{2}m_{p+1/2}(t)+ \frac{\gamma_{ p}}{2}\left(m_{1/2}(t)M_p +m_0(t) M_{p+1/2}\right) +\gamma_p \left(S_p(t)+\widetilde{S}_p(t)\right)
\end{multline*}
which is enough to prove that moments are uniformly propagated with time independently of $\alpha$, i.e., for any $p \geq 1$, there exists $C_p >0$ (independent of $\alpha$) such that $$m_p(0) < \infty \implies \sup_{t \geq 0}m_p(t) \leq C_p.$$
Let us now introduce the renormalized moments
$${z}_p(t):=\frac{ {m}_p(t)}{\Gamma(ap+b)}, \ \ \mbox{with}\ \ a=2/s,$$
where $\Gamma(\cdot)$ denotes the Gamma function and $b >0$ is a parameter to be fixed later on. Notice that, to get \eqref{Hadamard}, it suffices now to prove that, for any $a \geq 1$, one can find $b >0$ and  some positive constant $K >0$ (both $b$ and $K$ independent of $\alpha$) such that
\begin{equation}\label{Hadamard2} z_p(t) \leq K^p \qquad \forall p \geq 1.\end{equation}
 An important simplification, first observed in \cite{BoGaPa}, consists in noticing that, for any $a \geq 1$ and $b >0$,
\begin{equation*}
{S}_{p}(t)\leq C\;\Gamma(ap+a/2+2b)\;\mathcal{Z}_{p}(t)
\end{equation*}
where $C=C(a,b) >0$ does not depend on $p$ and
\begin{equation*}
\mathcal{Z}_p(t)=\max_{1\leq k\leq k_p}\left\{z_{k+1/2}(t)\;z_{p-k}(t),z_k(t)\;z_{p-k+1/2}(t)\right\}.
\end{equation*}
In the same way, one proves easily that
\begin{equation*}
\widetilde{S}_{p}(t)\leq C\;\Gamma(ap+a/2+2b)\;\widetilde{\mathcal{Z}}_{p}(t)\ \ \mbox{for}\ \ \ a\geq1,\;b>0,
\end{equation*}
where \begin{equation*}
\widetilde{\mathcal{Z}}_p(t)=\max_{1\leq k\leq k_p}\left\{z_{k+1/2}(t)\;\zeta_{p-k},\;z_{p-k+1/2}(t)\zeta_k \right\}
\end{equation*}
with $\zeta_p=\dfrac{M_p}{\Gamma(ap+b)}.$ Then, using the approximation formula
$$\lim_{p \to \infty}\dfrac{\Gamma(ap+r)}{\Gamma(ap+t)}(ap)^{t-r}=1  \qquad \forall a,t,r >0$$
together with the fact that $\gamma_{p}=O\left(\frac{1}{p}\right)$ as $p \to \infty$, we get that, for sufficiently large $p_* \geq 1$, there exist $c_1,c_2,c_3, c_4 >0$ such that
\begin{multline}\label{zpt}
\frac{\d z_p}{\d t}(t)+c_1\;p^{a/2}z_{p}(t)^{1+1/2p}\leq c_2\left(p^{-1}\zeta_p+p^{-1+a/2}\zeta_{p+1/2}\right) +c_3 p^{\tfrac{a}{2}+b-1}\widetilde{\mathcal{Z}}_p(t)\\
+c_4 \,p^{a/2+b-1}\;\mathcal{Z}_p(t)  \qquad \forall t \geq 0, \:p \geq p_*.
\end{multline}
Now, since $\M_0$  has an exponential tail of order $2$, \textit{a fortiori} it has an exponential tail of order $s$ with $0 < s \leq 2.$ Thus, for any $a=2/s \geq 1$ and  any $b >0$, there exists $C(a,b)>0$
and $A >1$ such that
$$M_p \leq C(a,b) \Gamma(ap+b) A ^p, \qquad \forall p \geq 1,$$
i.e. $\zeta_p \leq C(a,b) A ^p $  for any $p \geq 1.$ Therefore, \eqref{zpt} becomes
\begin{multline}\label{zptt}
\frac{\d z_p}{\d t}(t)+C_1\;p^{a/2}z_{p}(t)^{1+1/2p}\leq C_2\left(p^{-1}A^p+p^{-1+a/2}A^{p+1/2}\right) +C_3 p^{\tfrac{a}{2}+b-1}\mathbf{Z}_p(t)\\
+C_4 \,p^{\tfrac{a}{2}+b-1}\;\mathcal{Z}_p(t)  \qquad \forall t \geq 0, \:p \geq p_*
\end{multline}
for some positive constants $C_1,C_2,C_3,C_4 >0$
where $$\mathbf{Z}_p(t)=\max_{1\leq k\leq
  k_p}\left\{z_{k+1/2}(t)\;A^{p-k},\;z_{p-k+1/2}(t)A^k \right\}.$$ The
key observation is that, for any $p \geq p_*$, the functions
$\mathbf{Z}_p(t)$ and $\mathcal{Z}_p(t)$ involve $z_k(t)$ for $k \leq
p-1/2$ and \textit{do not involve} $z_p(t)$. This is the reason why
we will argue by \textit{induction} in order to prove that, for any $a
\geq 1$, if we choose $0<b<1$ it is possible to find $K>0$ large
enough so that $z_p(t)\leq K^{p}$. First, because of the exponential
integrability assumption on the initial datum $f_0,$ there exists $K_0
>0$ such that $z_p(0) \leq K_0^p$ for any $p \geq 1.$ Let us consider
now $p_0 \geq p_* > 1$ such that
\begin{equation*}
2C_2 p^{-1}_0 + (C_3+C_4) p_0^{b-1} \leq C_1,
\end{equation*}
and let $K >0$ be such that
\begin{equation*}
K\geq\left\{\max_{1\leq k\leq p_0}\sup_{ t\geq0} z_k(t),K_0, 1,A\right\}.
\end{equation*}
 Since moments of $f(t,v)$ are uniformly propagated, the existence of such a \textit{finite} $K$ is guaranteed. Defining now
    $$y_p(t):=K^{p} \qquad \forall t \geq 0$$ one can prove by induction (using also standard comparison of ODE's) that, for any $p \geq p_0$ with $2p \in \mathbb{N}$, $y_p(t)$  satisfies the differential inequality
\begin{equation*}\begin{split}
\frac{\d y_p}{\d t}(t)&+C_1\;p^{a/2}y_p(t)^{1+1/2p}\geq C_2\left(p^{-1}A^p+p^{-1+a/2}A^{p+1/2}\right) \\
& + C_3\;p^{\tfrac{a}{2}+b-1}\;\mathbf{Z}_p(t)+ C_4 \,p^{a/2+b-1}\;\mathcal{Z}_p(t),
\end{split}
\end{equation*} with moreover $y_p(0)\geq z_p(0)$. One deduces from this that $z_p(t) \leq y_p(t)=K^p$ for any $p\geq p_0$ and any $t \geq 0$. Notice  that the comparison argument for ODE's is licit here since, again, for a given $p \geq p_0$,  $\mathbf{Z}_p(t)$ and $\mathcal{Z}_p(t)$ involve only $z_k(t)$ for $k \leq p-1/2$. This yields the desired conclusion \eqref{Hadamard}.\end{proof}

We can now give a stationary version of the above Theorem in order to deduce the order of the exponential tail of the solutions to \eqref{equi}. In the elastic case $\alpha=1$, as we already saw it, the solution to \eqref{equi} is a Maxwellian and therefore has an exponential tail of order $2$.
For a given $\alpha \in (0,1)$, we look for the order of the exponential tail of the solution $F_\alpha$ to \eqref{equi}. Notice that, as shown in \cite{MiMo}, the bounds obtained from $\Q_\alpha$ are actually uniform with respect to the coefficient $\alpha$.
This suggests that the order of the exponential tail of the solution $F_\alpha$ shall be independent of $\alpha$.
Since, for $\alpha=1$, the order is $s=2$, we infer that the solution $F_\alpha$ to \eqref{equi} has an {\it exponential
tail of order} $2.$ This is the object of the following Theorem:
\begin{theo}\label{gaus}  There exist  some constant $A >0$ and $M >0$ such that, for any  $\alpha \in (0,1]$ and any solution $F_\alpha$ to \eqref{equi} one has
$$\int_{\R^3} F_\alpha(v) \exp\left(A |v|^2\right) \d\v \leq M.$$
\end{theo}
\begin{proof} The proof of Theorem \ref{gaus} follows exactly the same lines as that of Theorem \ref{tailT}. We only sketch here the straightforward modifications. Recall that, for any $r,s >0$ and any $\alpha \in (0,1]$, we defined
\begin{equation}
\label{series_rs}
\mathscr{F}_{r,s}(\fe) = \int_{\R^3} \fe \,
\Big( \sum\limits_{k=0}^\infty
\,\frac{r^k}{k!} \,|v|^{sk}\Big)\, \d v
=  \sum_{k=0}^\infty \,
\frac{r^k }{k!}\, \,\mathbf{m}_{\frac{sk}{2}}(\alpha)
\end{equation}
where
\begin{equation}
\label{moments}
\mathbf{m}_{p}(\alpha)=\IR F_\alpha(v) |v|^{2p}\, dv, \quad p\geq 0.
\end{equation}
For any $p \geq 0$, we introduce now the following stationary moments
$$Q_p(\alpha)=\IR \Q_\alpha(F_\alpha,F_\alpha)(v)|v|^{2p}\d\v,\qquad \qquad L_p(\alpha)=\IR \L(F_\alpha)|v|^{2p}\d\v.$$
Of course, for any $p \geq 0,$ $Q_p(\alpha)+L_p(\alpha)=0.$  Arguing exactly as above we get that
\begin{equation}\label{estimate1}
3(1-\gamma_p) \mathbf{m}_{p}^{1+1/2p}(\alpha) \leq {\gamma_p}  \left(S_p(\alpha)+\widetilde{S}_p(\alpha) + \mathbf{m}_{1/2}(\alpha)M_p +  M_{p+1/2}\right)
\end{equation}
where
$$S_p(\alpha)=\sum_{k=1}^{k_p} \binom{p}{k}
\left(\mathbf{m}_{ k+1/2 }(\alpha)\mathbf{m}_{ p-k }(\alpha)+\mathbf{m}_{ p-k+1/2 }(\alpha) \mathbf{m}_{ k}(\alpha)\right)$$
while
\begin{equation*}
  \widetilde{S}_p(\alpha)=\sum_{k=1}^{k_p}
  \binom{p}{k} \left(
    \mathbf{m}_{ k+1/2 } (\alpha) M_{ p-k }
    +\mathbf{m}_{ p-k+1/2 }(\alpha) M_{ k}
  \right).
\end{equation*}
To prove that the solution $F_\alpha$ to \eqref{equi} has an {\it
  exponential tail of order} $2$, as in the above proof
(with $s=2$) it is sufficient to prove that there exist $C >0$ and $X
>0$ such that
\begin{equation}\label{gammap}
\mathbf{m}_p(\alpha) \leq C \Gamma(p+1/2) X^p, \qquad \forall p \geq 1,\quad \forall \alpha \in (0,1].\end{equation}
Notice that, since obviously $\M_0$ has an exponential tail of order $2$, there exists $C_0>0$
and $X_0 >1$ such that
$$M_p \leq C_0 \Gamma(p+1/2) X_0^p, \qquad \forall p \geq 1.$$
Then, arguing as in the above proof, one gets that the above decrease of $M_p$ is enough to get \eqref{gammap} by an induction argument as in the proof of Theorem \ref{tailT}.\end{proof}

\begin{nb} The above Theorem provides some weighted $L^1$ space which contains
all the stationary solutions for any $\alpha \in (0,1)$. Moreover, since the conclusion of the above result should hold for any $\alpha \in (0,1)$, in particular, for $\alpha=1$ since
$F_\alpha=\M_1$ is an explicit Maxwellian, one has $A < A^\sharp$ with $A^\sharp:=\dfrac{1}{2\Theta^\sharp}.$
\end{nb}

\subsection{Uniform bound for the energy and control of the $L^2$-norm}
Upper bounds for the energy of the solution to \eqref{equi} are easily obtained as a consequence of the above calculations or, more simply,
from \cite[Eq. (4.6)]{BiCaLo}: there exists $E_\mathrm{max} >0$ such that
$$E_\alpha:=\IR |v|^2 F_\alpha(v)\d\v < E_\mathrm{max} \qquad \qquad \forall \alpha \in (0,1]$$
where $F_\alpha$ is a solution to \eqref{equi}. In order to derive an uniform lower bound of $E_\alpha$
(showing in particular that $E_\alpha$ does not vanish in the elastic limit $\alpha \to 1$), one shall actually derive
an uniform lower bound of the $L^2$-norm of any solution to \eqref{equi}:
\begin{theo}\label{theoenergy} Given $\alpha \in (0,1]$, any stationary solution $F_\alpha$ to \eqref{equi} belongs to $L^2(\R^3,\d\v).$
More precisely, there exists an uniform constant $\ell_2>0$ such that
$$\|F_\alpha\|_{L^2(\R^3,\d\v)} \leq \ell_2.$$
As a consequence, there exists $E_\mathrm{min} >0$ such that
$$E_\mathrm{min} \leq E_\alpha \leq E_\mathrm{max} \qquad\qquad \forall \alpha \in (0,1].$$
\end{theo}
\begin{proof} We prove the control of the $L^2$ norm as in \cite{MiMo2}. Precisely, let $A >0$ be fixed and let
$\Lambda_A(x)=\frac{x^2}{2}\chi_{\left\{x <A\right\}}+(Ax-\frac{A^2}{2})\chi_{\left\{x>A\right\}},$ $x \in \R.$
The function $\Lambda_A$ is a $\mathcal{C}^1$-function over~$\R$ and $\lim_{A\to \infty}\Lambda_A(x)=\frac{x^2}{2}$ for any $x \in \R.$
In particular, for proving the claim, it is enough proving that there exists some positive constant $c >0$ not depending on
$\alpha \in (0,1]$ such that
\begin{equation}\label{LAMA}\limsup_{A \to \infty} \IR \Lambda_A\left(F_\alpha\right)(v)\d\v \leq c.\end{equation}
Let $T_A(x):=\min\left(x,A\right)=\Lambda'_A(x)$. Multiplying the identity \eqref{equi} by $T_A(F_\alpha)$ and integrating over $\R^3$ leads to
\begin{multline*}
\IR \Q_\alpha^-(F_\alpha,F_\alpha)T_A(F_\alpha)\d\v + \IR \L^-(F_\alpha)T_A(F_\alpha)\d\v\\
= \IR \Q_\alpha^+(F_\alpha,F_\alpha)T_A(F_\alpha)\d\v + \IR \L^+(F_\alpha)T_A(F_\alpha)\d\v \qquad \alpha \in (0,1].
\end{multline*}
All the integrals in the above expression are nonnegative and in particular:
$$ \IR \L^-(F_\alpha)T_A(F_\alpha)\d\v\leq \IR \Q_\alpha^+(F_\alpha,F_\alpha)T_A(F_\alpha)\d\v + \IR \L^+(F_\alpha)
T_A(F_\alpha)\d\v.
$$
Now, one estimates the left-hand side from below \textit{uniformly with respect to} $\alpha$ as follows:
\begin{multline*}\IR \L^-(F_\alpha)T_A(F_\alpha)\d\v=\IR F_\alpha(v)T_A(F_\alpha)(v)\d\v\IR \M_0(w)|v-w|\d\w
\\\geq c_M \IR F_\alpha(v)T_A(F_\alpha)(v)(1+|v|)\d\v\end{multline*}
where $c_M=\inf_v \IR \M_0(w)|v-w|\d\w/ (1+|v|)$ is positive and finite. In particular, it does not depend on $\alpha.$
Then, as in \cite{MiMo2}, since $\Lambda_A(x) \leq x T_A(x)$ one gets that
$$c_M\IR \Lambda_A(F_\alpha)(1+|v|)\d\v\leq \IR \Q_\alpha^+(F_\alpha,F_\alpha)T_A(F_\alpha)\d\v +
\IR \L^+(F_\alpha)T_A(F_\alpha)\d\v.
$$
Now, according to \cite[Step 2, Proposition 2.1]{MiMo2}, there exists $\theta \in (0,1)$ such that, for any $\alpha \in (0,1]$, there is a constant
$\mathcal{C}=\mathcal{C}(E_\alpha) >0$ and $A_\alpha >0$ such that
$$\IR \Q_\alpha^+(F_\alpha,F_\alpha)T_A(F_\alpha)\d\v \leq \mathcal{C} \|T_A(F_\alpha)\|_{L^2}^{2(1-\theta)} +
\dfrac{c_M}{2} \IR \Lambda_A(F_\alpha)(1+|v|)\d\v \qquad \forall A >A_\alpha.$$
Notice that, though $A_\alpha$ depends on the inelasticity parameter $\alpha$, it will play no role since we are only considering
the limit as $A$ goes to infinity.
Moreover, a careful reading of the proof of \cite[Prop. 2.1]{MiMo2} shows that the constant
$\mathcal{C}$ depends on $\alpha$ only through upper bounds of the energy $E_\alpha$. In particular,
since we proved that $E_\alpha \leq E_\mathrm{max}$,
one can set $\mathcal{C}=\sup_{0 < \alpha < 1} \mathcal{C}(E_\alpha) < \infty$
in the above inequality. One obtains finally
$$\dfrac{c_M}{2} \IR \Lambda_A(F_\alpha)(1+|v|)\d\v \leq \mathcal{C} \|T_A(F_\alpha)\|_{L^2}^{2(1-\theta)}
+ \IR \L^+(F_\alpha)T_A(F_\alpha)\d\v \qquad \forall A > A_\alpha.$$
One estimates now the last integral on the right-hand side owing to
\cite[Theorem 1]{AlCaGa}. Precisely, according to the Cauchy-Schwarz
inequality
$$\IR \L^+(F_\alpha)T_A(F_\alpha)\d\v \leq \|\L^+(F_\alpha)\|_{L^2} \|T_A(F_\alpha)\|_{L^2} \qquad \forall A >0$$
(notice that $T_A(F_\alpha) \in L^2$ for any fixed $A$ since $\big[ T_A(x) \big]^2 \leq Ax$). Now, using the fact that $\L^+(f)=\Q_e^+(f,\M_0)$,
one deduces directly from \cite[Theorem 1]{AlCaGa} that
$$\|\L^+(F_\alpha)\|_{L^2} \leq C_e \|F_\alpha\|_{L^1_1} \|\M_0\|_{L^2_1}$$
where $C_e >0$ depends on the inelasticity parameter $e$. Since $\sup_{0 <\alpha <1}\|F_\alpha\|_{L^1_1} < \infty$
according to the result of the previous section, one gets that
$$\IR \L^+(F_\alpha)T_A(F_\alpha)\d\v \leq  C_0 \|T_A(F_\alpha)\|_{L^2} \qquad \forall A >0$$
where $C_0 >0$ is a positive constant independent of $\alpha$. Finally, we obtain
$$\dfrac{c_M}{2} \IR \Lambda_A(F_\alpha)(1+|v|)\d\v \leq \mathcal{C} \|T_A(F_\alpha)\|_{L^2}^{2(1-\theta)}
+ C_0 \|T_A(F_\alpha)\|_{L^2} \qquad \forall A >A_\alpha.$$
Since $\Lambda_A(F_\alpha) \geq T_A(F_\alpha)^2/2$, this means that
$$\dfrac{c_M}{4} \|T_A(F_\alpha)\|_{L^2}^2 \leq \mathcal{C} \|T_A(F_\alpha)\|_{L^2}^{2(1-\theta)}
+ C_0 \|T_A(F_\alpha)\|_{L^2} \qquad \forall A >A_\alpha$$
which clearly implies \eqref{LAMA}. Now, it is a classical feature to deduce the uniform lower bound of the energy $E_\alpha$ from the uniform control
of $\|F_\alpha\|_{L^2}$ (see, e.g. \cite[Proposition 4.6]{BiCaLo}).
\end{proof}
As in \cite[Prop. 2.1, step 8]{MiMo2}, a simple corollary of the above Theorem and the results of \cite[Section 6, Prop. 6.2]{BiCaLo} is the following uniform smoothness estimate:
\begin{cor}\label{sobolev} For any $k \in \mathbb{N}$, there exists $C_k >0$ such that
$$\|\fe\|_{H^k(\R^3)} \leq C_k, \qquad \forall \alpha \in (0,1].$$
In particular, there exists $C_\infty >0$ such that
$$\|\fe\|_{L^\infty(\R^3)} \leq C_\infty \qquad \forall \alpha \in (0,1].$$
\end{cor}
\begin{proof} Notice that, from the uniform lower bound on the energy $E_\alpha \geq E_\mathrm{min}$ $(\alpha \in (0,1])$, one notices that there exists some positive constant $c_0 >0$ such that
$$\Sigma(\fe)(v)=\int_{\R^3} |\v-\w|\fe(\w)\d\w \geq c_0 \left(1+|v|\right) \qquad \forall \alpha \in (0,1].$$
This, together with \cite[Proposition 6.2]{BiCaLo},  is enough to provide the uniform bound in $H^k(\R^3)$ as in \cite[Prop. 2.1, step 8]{MiMo2}.
Now, by Sobolev embedding theorem, one gets the second part of the corollary choosing simply $k>6$.
\end{proof}

\section{Uniform pointwise estimates}

On the basis of the previous result, we derive in this section pointwise estimates for the steady state $\fe$ which are uniform with respect to the inelasticity parameter $\alpha.$ More precisely, we shall prove that there exists two Maxwellian distributions $\underline{\M}$ and $\overline{\M}$ (independent of $\alpha$) such that
$$\underline{\M}(v) \leq \fe(v) \leq \overline{\M}(v) \qquad \forall v \in \R^3.$$
We will treat separately the upper bound and the lower bound.
\subsection{Uniform pointwise upper Maxwellian bound}  The strategy of proof is inspired by the comparison principle of \cite{GaPaViM} and uses some estimates of \cite{AlCaGa}. Precisely, a first general comparison principle is the following:
\begin{propo}\label{maxi} Let $\alpha \in (0,1]$ and $\fe$ be a solution to \eqref{equi} with unit mass. Assume there exists a measurable subset $\mathcal{U}$ of $\R^3$ (with nonzero Lebesgue measure) and a measurable and nonnegative distribution $G=G(v)$ such that
\begin{equation}\label{U<}
\Q_\alpha(G,\fe) + \L(G) < 0 \qquad \text{ for any } \quad v \in \mathcal{U};\end{equation}
and
\begin{equation}\label{U>}
\fe(v) \leq G(v) \qquad \text{ for any } \quad v \in \R^3 \setminus \mathcal{U}.
\end{equation}
Then, $\fe(v) \leq G(v)$ for almost every $v \in \R^3.$
\end{propo}
\begin{proof} As already said, the proof follows the strategy of
  \cite[Theorem 3]{GaPaViM} which is given for the time-dependent
  (space inhomogeneous) elastic Boltzmann equation. We adapt it in a
  simple way for granular gases in spatially homogeneous
  situations. One notices first that, for any nonnegative distribution
  $g \geq 0$ and any distribution $f$
\begin{equation}\label{negaqa}
\IR \Q_\alpha(f,g)(v)\mathrm{sign}(f)(v) \d v \leq 0.\end{equation}
Indeed, according to \eqref{weakfg}, the above integral is equal to
$$\frac{1}{4\pi}\IRR\IS  f(v)g(\vb)|v-\vb|\left(\mathrm{sign}(f)(v')-\mathrm{sign}(f(v))\right)\d v\d\vb\d\sigma$$
and the conclusion follows since $g(\vb) \geq 0$ while $f(v)\left(\mathrm{sign}(f)(v')-\mathrm{sign}(f(v))\right) \leq 0$ for any $v,\vb \in \R^3.$ For the same reason
\begin{equation}\label{negal}
\IR \L(f)(v)\mathrm{sign}(f(v)) \d v \leq 0 \qquad \text{ for any distribution } f.
\end{equation}
Now,  after multiplying \eqref{equi} by $\mathrm{sign}(\fe-G)$ and integrating over $\R^3$, one gets
\begin{equation*}\begin{split}
0&=\IR \left(\Q_\alpha(\fe,\fe)+\L(\fe)\right)\,\mathrm{sign}(\fe-G)\d v=\IR \left(\Q_\alpha(G,\fe)+\L(G)\right)\,\mathrm{sign}(\fe-G)\d v\\
 &\phantom{++++} + \IR \left(\Q_\alpha(\fe-G,\fe)+\L(\fe-G)\right)\,\mathrm{sign}(\fe-G)\d v
\end{split}
\end{equation*}
and this last integral is nonpositive according to \eqref{negaqa} and \eqref{negal}. Therefore,
$$\IR  \left(\Q_\alpha(G,\fe)+\L(G)\right)\,\mathrm{sign}(\fe-G)\d v \geq 0.$$
Moreover, using the  mass conservation property of both the collision operators, one can rewrite the above inequality as
$$\IR  \left(\Q_\alpha(G,\fe)+\L(G)\right) \,\dfrac{\mathrm{sign}(\fe-G)+1}{2}\d v \geq 0$$
where $(\mathrm{sign}(\fe-G)+1)/2$ is always nonnegative. We split this integral over $\mathcal{U}$ and its complementary. Whenever $v \notin \mathcal{U},$ by assumption \eqref{U>} one has $(\mathrm{sign}(\fe(v)-G(v))+1)/2=0.$ Thus, the integral over $\R^3$ reduces to the integral over $\mathcal{U}$, i.e.
$$\int_{\mathcal{U}} \left(\Q_\alpha(G,\fe)+\L(G)\right) \,\dfrac{\mathrm{sign}(\fe-G)+1}{2}\d v \geq 0.$$
Now, according to our assumption \eqref{U<}, the above is the integral of a nonpositive measurable distribution. Therefore,
$$\left(\Q_\alpha(G,\fe)+\L(G)\right) \,\left(\mathrm{sign}(\fe-G)+1\right)=0 \qquad \text{ almost everywhere over } \mathcal{U}.$$
Using again \eqref{U<} we get that $\mathrm{sign}(\fe(v)-G(v)) =-1$ for almost every $v \in \mathcal{U}$ which proves that $\fe(v) \leq G(v)$ for almost every $v \in \R^3.$
\end{proof}
Now, in order to prove that every steady state $\fe$ is bounded from above  by an universal Maxwellian distribution, we only have to determine a Maxwellian distribution $G$ and a measurable subset $\mathcal{U}$ for which the above \eqref{U<} and \eqref{U>} hold true. We will need the following general result, proven in \cite[Proposition 11]{AlCaGa} that we state here for hard-spheres interactions only:
\begin{theo}[\textbf{\textit{Alonso et al.}\cite{AlCaGa}}] Let $1 \leq  p, q, r \leq \infty$ with $1/p + 1/q = 1 + 1/r.$ Then, for $a > 0$ there is a positive constant $C_a >0$ such that
\begin{equation}\label{Alo}
\left\|\Q_\alpha^+(f,g)\sM_a^{-1}\right\|_{L^r(\R^3)} \leq C_a \left\|f \sM_a^{-1}\right\|_{L^p(\R^3)}\,\left\|g\sM_a^{-1}\right\|_{L^q_1(\R^3)} \qquad \forall \alpha \in (0,1]\end{equation}
where  $\sM_a(v)=\exp(-a|v|^2)$, $v \in \R^3.$
\end{theo}
\begin{nb} Notice that in \cite{AlCaGa} the constant appearing in
  \eqref{Alo} is actually given by $C_\alpha C_{1,a}$, where $C_{1,a}$
  is given by \cite[Eq. (6.10)]{AlCaGa} and depends only on $a$, while
  $C_\alpha$ is given by \cite[Eq.~(4.4)]{AlCaGa} and depends on the
  inelasticity parameter only through $(1-\alpha)^2.$ Bounding this
  last quantity simply by $1$, one sees that $\sup_{\alpha \in
    (0,1]}C_\alpha < \infty$, thus obtaining a constant $C_a$ in~\eqref{Alo} which does not depend on the inelasticity parameter
  $\alpha.$
\end{nb}
This leads to the following
\begin{theo} For any positive number $a < \min(\frac{1}{2\Theta_0},A)$ where $A >0$ is given in Theorem~\ref{gaus}, there exists a positive constant $\mu_a >0$ (independent of the inelasticity parameter $\alpha$) such that
$$\fe(v) \leq \exp\big(-a|v|^2 + \mu_a \big)  \qquad \forall v \in \R^3, \forall \alpha \in (0,1].$$
\end{theo}
\begin{proof} Let us fix $a<\min(\frac{1}{2\Theta_0},A) $ and set $\sM_a(v)=\exp(-a|v|^2)$, $v \in \R^3.$
As in \cite{GaPaViM}, one shall apply Proposition \ref{maxi} with
$$\mathcal{U}=\{v \in \R^3\,,\,|v| >R\}$$
for $R >0$ sufficiently large and with $G(v)=K_a \sM_a(v)$ and $K_a$ to be determined. The technical part is to prove that \eqref{U<} holds true for $R >0$ large enough. First, one has
$$\Q_\alpha^-(\sM_a,\fe)(v)+\L^-(\sM_a)(v)=\sM_a(v)\left(\IR \fe(w)|v-w|\d w + \IR \M_0(w)|v-w|\d w\right)$$
Recall that according to Theorem \ref{gaus}, $\sup_{\alpha \in (0,1]} \IR \fe(v)|v|\d v=\overline{m}_1 <\infty.$ Therefore, since $|v-w| > |v|-|w|$ and both $\fe$ and $\M_0$ have unit mass, one has in a direct way
\begin{equation}\label{Q-Ma}\Q_\alpha^-(\sM_a,\fe)(v)+\L^-(\sM_a)(v) \geq 2\sM_a(v) \left(|v|-\frac{\overline{m}_1 +\overline{m}_0}{2}\right) \qquad \forall v \in \R^3\end{equation}
where $\overline{m}_0=\IR \M_0(w)|w|\d w.$ Now, to estimate $\L^+(\sM_a)=\Q^+_e(\sM_a,\M_0)$, one applies  \eqref{Alo} with  $f=\sM_a$, $g=\M_0$ and $(p,q,r)=(\infty,1,\infty).$ Since $a < \frac{1}{2\Theta_0}$, one sees that $\|\M_0 \sM_a^{-1}\|_{L^1_1(\R^3)} < \infty$ while trivially $\|\sM_a \sM_a^{-1}\|_{L^\infty(\R^3)}=1$. Thus, there exists a positive constant $c_1(a) >0$ such that
$$\L^+(\sM_a)(v) \leq c_1(a)\sM_a(v) \qquad \forall v \in \R^3.$$
To estimate $\Q_\alpha^+(\sM_a,\fe)$, one applies now \eqref{Alo} with $f=\sM_a$, $g=\fe$ and $(p,q,r)=(\infty,1,\infty).$ Since
$\sup_{\alpha} \IR \fe(v)\exp(A|v|^2)\d v <\infty$
according to Theorem \ref{gaus}, one sees that, for any $a < A$,
$$\sup_{\alpha \in (0,1]} \|\fe \sM_a^{-1}\|_{L^1_1(\R^3)} < \infty$$
and therefore there is a positive constant $c_2(a) >0$ such that $\Q_\alpha^+(\sM_a,\fe)(v) \leq c_2(a)\sM_a(v)$ $\forall v \in \R^3.$ Gathering these two estimates, one gets the existence of a positive constant $C_a$ (independent of $\alpha$) such that
\begin{equation}\label{Q+Ma}
\Q_\alpha^+(\sM_a,\fe)(v)+\L^+(\sM_a)(v) \leq C_a \sM_a(v) \qquad \forall v \in \R^3.
\end{equation}
Combining \eqref{Q-Ma} and \eqref{Q+Ma}, one sees that, choosing $R > \frac{\overline{m}_1+\overline{m}_0+C_a}{2}$, we have
$$\Q_\alpha^-(\sM_a,\fe)(v)+\L^-(\sM_a)(v) \geq C_a \sM_a(v) \geq \Q_\alpha^+(\sM_a,\fe)(v) + \L^+(\sM_a)(v) \qquad \forall |v| > R
$$
i.e.
$$\Q_\alpha(\sM_a,\fe)(v)+ \L(\sM_a) (v) \leq 0 \qquad \forall v \in \mathcal{U}.$$
Now, since there exists $C >0$ such that $\fe(v) \leq C$ for any $v \in \R^3$ and any $\alpha \in (0,1]$ according to Corollary \ref{sobolev}, it is clear that one can find a positive constant $K_a=C\exp(-aR^2) >0$ such that
$$\fe(v) \leq K_a \sM_a(v) \qquad \forall |v| \leq R.$$
With this choice of $R$ and $K_a$, the function $G=K_a \sM_a(v)$ satisfies \eqref{U<} and \eqref{U>} of Proposition \ref{maxi} and we get our conclusion with $\mu_a=\log K_a.$\end{proof}

\subsection{Uniform pointwise lower Maxwellian bound}

We prove now a Maxwellian pointwise lower bound for the stationary
solution $\fe$ which is uniform with respect to the inelasticity
parameter. It turns out that the proof of such a result is much
simpler than the ones yielding (non Maxwellian) pointwise lower bounds
in the diffusively driven case~\cite{MiMo3} or for the homogeneous
cooling state in \cite{MiMo2}. These two results rely on the spreading
properties of the \textit{nonlinear} inelastic collision operator
$\Q_\alpha$ (in the spirit of similar results obtained in the elastic
case in \cite{ada}). On the contrary, our approach relies uniquely on
the properties of the \textit{linear} collision operator $\L$ and,
more precisely, on the explicit integral representation of $\L^+$
derived in \cite{ArLo}. We first prove a general lower bound for the
time-dependent problem \eqref{BEev}:
\begin{theo}\label{inf2t} Let $f_0 \in L^1_3$ be a nonnegative initial datum with unit mass and let $f(t,v)$ be the associated solution to \eqref{BEev}.
 Then, for any $t_0 >0$, there exists a positive constant $a_0 >0$ (which depends only on $C >0$ and $t_0$) such that
$$f(t,v) \geq a_0  \exp(-\gamma_1 |v|^2) \qquad \forall v \in \R^3\,,\forall t \geq t_0$$
where $\gamma_1=\frac{3+3\mu+\mu^2}{4\Theta_0}$ and $\mu=2 \frac{1-e}{1+e}$.
\end{theo}
\begin{proof} The solution $f(t,v)$ to \eqref{BEev} satisfies
$$\partial_t f(t,v) + \left(\Sigma(f(t))(v)+\sigma(v)\right)f(t,v) =\Q^+_\alpha(f(t)\,,\,f(t))(v)+\L^+(f(t,\cdot))(v).$$
Moreover, because of the propagation of moments uniformly with respect to $\alpha$, there is some $M_1 >0$, independent of $\alpha$ such that
\begin{equation}\label{M1} \sup_{t\geq 0} \IR f(t,v)|v|\d v \leq M_1 < \infty\end{equation}
 so that  there exists $c_2 >0$ such that
$$ \Sigma(f(t))(v) + \sigma(v) \leq c_2(1+|v|) \qquad \forall  v \in \R^3, \:\forall t \geq 0.$$
Therefore, the solution $f(t,v)$ satisfies the following inequality:
\begin{equation}\label{inequalityft}\partial_t f(t,v) + c_2(1+|v|)f(t,v) \geq \L^+(f)(t,v) \qquad \forall t \geq 0\,,\:v \in \R^3.\end{equation}
Now, according to \cite{ArLo}, the positive part $\L^+$ admits the following integral representation
$$\L^+(f)(t,v)=\IR k(v,w)f(t,w)\d w$$
where \begin{equation}\label{k}
k( v,w)=C_0|v-w|^{-1}\exp
\left\{-\beta_0 \left(
(1+\mu)|v-w|+\dfrac{|v|^2-{|w|}^2}{|v-w|}\right)^2\right\}\end{equation}
with $\mu=2 \frac{1-e}{1+e} \geq 0$, $\beta_0=\frac{1}{8\Theta_0}$ and $C_0 >0$ is a positive constant (depending on $e$ and $\Theta_0$). Moreover,  the \textit{microscopic detailed balance law} holds true
$$k(v,w)\M (w)=k(w,v)\M(v) \qquad \forall v,w \in \R^3$$
where $\M(v)$ is the Maxwellian distribution defined in \eqref{steadystate}:
$$\M(v)=\left(\dfrac{1}{2\pi
\Theta^{\#}}\right)^{3/2}\exp(-A^\sharp |v|^2), \quad A^\sharp=\frac{1}{2\Theta^\sharp}=4(1+\mu)\beta_0.$$
Therefore $\L^+(f)(t,v)=\M(v)\ds\IR k(w,v)f(t,w)\M^{-1}(w)\d w.$ Since
$$|v-w|^2 \leq 2|v|^2+2|w|^2 \qquad \text{ and }  \qquad\dfrac{(|w|^2-|v|^2)^2}{|v-w|^2} \leq 2|v|^2+2|w|^2$$
straightforward computations yield
$$k(w,v) \geq \dfrac{C_0}{|v|+|w|}\,\exp(-\gamma_0|v|^2)\;\exp(-\gamma_1|w|^2),$$
with $\gamma_0=2\beta_0(1+\mu+\mu^2)$ and $\gamma_1=2\beta_0(3+3\mu+\mu^2)$. Notice that  $A^\sharp-\gamma_1=-\gamma_0.$  Now, owing to the mass condition and \eqref{M1}, for any $R  \geq 2M_1$ one has
\begin{equation}\label{BORt}\inf_{t \geq 0} \int_{B(0,R)}f(t,w)\d w \geq {1 \over 2}\end{equation}
and
\begin{equation*}\begin{split}
\L^+(f)(t,v) &\geq \dfrac{C_0}{(2\pi
\Theta^{\#})^{3/2}} \exp(-(A^\sharp+\gamma_0)|v|^2)\times \\
&\phantom{++++}\int_{B(0,R)}\exp(-\gamma_1|w|^2)\M^{-1}(w)f(t,w) \dfrac{\d w}{|v|+|w|}\\
&= C_0 \exp(-(A^\sharp+\gamma_0)|v|^2)\times \\
&\phantom{++++}\int_{B(0,R)}\exp((A^\sharp-\gamma_1)|w|^2)f(t,w) \dfrac{\d w}{|v|+|w|}.
\end{split}\end{equation*}
Hence, there exists $C_R=C_0\exp(-\gamma_0 R{\,}^2) >0$ independent of $t \geq 0$ such that
$$\L^+(f)(t,v) \geq \dfrac{C_R}{|v|+R} \exp(-(A^\sharp+\gamma_0)|v|^2)\int_{B(0,R)} f(t,w)\d w.$$
This, together with \eqref{inequalityft} and \eqref{BORt}, yields
$$\partial_t f(t,v)  + c_2(1+|v|)f(t,v) \geq  \dfrac{C_R}{2(|v|+R)}\exp(-\gamma_1|v|^2) \qquad \forall v \in \R^3$$
from which we deduce
$$f(t,v) \geq  \dfrac{C_R}{2c_2(|v|+R)^2}\exp(-\gamma_1|v|^2)\left(1-e^{-c_2(1+|v|)t}\right) + e^{-c_2(1+|v|)t} f_0(v) \qquad \forall t \geq 0.$$
This clearly leads to the desired result.
\end{proof}
\begin{nb} Notice that the above proof does not require the energy of
  $f(t,v)$ to be bounded from below and the various constants involved
  depend only on the uniform upper bound on the first order moment
  \eqref{M1}. In particular, the lower bound of the previous Theorem
  \ref{inf2t} shows that there exists $a_1 >0$, independent of $\alpha
  \in (0,1]$ such that
$$\inf_{t \geq 0}\IR f_\alpha(t,v)|v|^2\d v \geq a_1 >0$$
for any solution $f_\alpha(t,v)$ to \eqref{BEev}.
\end{nb}

A stationary version of the above result is now straightforward:
\begin{theo}
  \label{inf2}
  There exists some positive constant  $a_0>0$ such that, for any $\alpha
  \in (0,1]$,
$$\fe(v) \geq a_0^{-1} \exp(-a_0|v|^2) \qquad \forall v \in \R^3.$$
\end{theo}
\begin{proof} The proof follows the same paths of the previous one and is omitted here. Notice that the constant $a_0 >0$ does not depend on $\alpha \in (0,1)$
because the bounds provided by  Theorem \ref{theoenergy} are uniform with respect to the inelasticity parameter.
\end{proof}

The above uniform lower bound together with the regularity estimates of Corollary~\ref{sobolev} has important consequences on the entropy production. Precisely, for any $\alpha \in (0,1]$, define the entropy dissipation functional, for any nonnegative~$g$:
\begin{multline*}\D_{H,\alpha}(g)=\frac{1}{8\pi} \int_{\R^3 \times \R^3\times  \mathbb{S}^2} |v-w| g(v)g(w)\\
\times\left(\dfrac{g(v')g(w')}{g(v)g(w)} - \log\dfrac{g(v')g(w')}{g(v)g(w)} -1\right)\d\sigma\d v \d w  \geq 0\end{multline*}
where the post-collisional velocities $(v',w')=(v_\alpha',w_\alpha')$ are defined in \eqref{co:transf}.
Notice that, for any nonnegative $g$ for which all the integrals make sense, one has  (see \cite{GPV**,MiMo2} for details):
\begin{equation} \label{dissDHal}
\IR \Q_\alpha(g,g)(v) \log g(v)\d v=-\D_{H,\alpha}(g) + \frac{1-\alpha^2}{\alpha^2}\IRR g(v)g(w)|v-w|\d v \d w
\end{equation}

Then, arguing exactly as in \cite[Corollary 3.4]{MiMo2}, we get the following
\begin{propo}\label{estimentropy} There exist $k_0$  and $q_0 \in \mathbb{N}$ large enough such that, for any $a_i >0$, there is some constant $C >0$ such that, for any $g$ satisfying
$$\|g\|_{H^{k_0} \cap L^1_{q_0}} \leq a_1, \qquad g(v) \geq a_2 \exp(-a_3 |v|^2)$$
one has
$$|\D_{H,\alpha}(g) -\D_{H,1}(g)| \leq C(1-\alpha) \qquad \forall \alpha \in (0,1].$$
\end{propo}
\begin{nb} Notice that, if $f_0 \in L^1_{q_0} \cap H^{k_0}$ is an initial distribution with unit mass, then, according to \cite[Proposition 6.3]{BiCaLo}, the associated solution $f(t,v)$ to \eqref{BEev} satisfies
$$\sup_{t \geq 0} \|f(t)\|_{H^{k_0} \cap L^1_{q_0}} \leq a_1$$
for some positive constant $a_1 >0$. Hence,  one deduces from Theorem \ref{inf2t} and the above Proposition that there exists some constant $C >0$ such that
$$\sup_{t \geq 0}|\D_{H,\alpha}(f(t)) -\D_{H,1}(f(t))| \leq C(1-\alpha) \qquad \forall \alpha \in (0,1].$$
\end{nb}
\section{Uniqueness of the steady state}
\label{sec:uniqueness}

We aim now to prove that, for $\alpha \in (0,1)$ large enough,
the steady state $F_\alpha$ is unique, precisely, we show there
is $\alpha_0 \in (0,1)$ such that, for any $\varrho > 0$ and
any $\alpha \in (\alpha_0,1)$, the set
\begin{equation}\label{Salpha}\mathfrak{S}_\alpha(\varrho)=\left\{F_\alpha \in L^1_2, \:F_\alpha \geq 0,\; \fe  \text{ solution to \eqref{equi} with } \IR \fe \d v =\varrho  \right\}\end{equation}
reduces to a singleton. We first recall that Theorem \ref{steadystate} proves that it is the case in the elastic case: $\mathfrak{S}_1(\varrho)=\{\varrho\M\}$ for any $\varrho >0.$ On the basis of this easy result,  we adopt the strategy described in the Introduction to prove the uniqueness of the steady state whenever $\alpha < 1.$   We begin by recalling the fundamental estimates of \cite{MiMo3} ensuring \eqref{XY} and \eqref{XY1}.
\subsection{Estimates on the collision operator} We recall here some results established in \cite{MiMo3} determining the function space in which the collision operator $\Q_\alpha$ depends continuously on the restitution coefficient $\alpha \in (0,1].$ Let
$$\mathcal{X} =L^1(m ^{-1})=L^1(\R^3,m^{-1}(v)\d v), \qquad \mathcal{Y} =L^1_1(m ^{-1})=L^1(\R^3, \langle v \rangle m^{-1}(v)\d v)$$
where
$$m(v)=\exp\left(-a |v|^s\right), \qquad a >0, s \in (0,1].$$
Then, from \cite[Proposition 3.2]{MiMo3}:
\begin{propo}[\textbf{Mischler \& Mouhot}]\label{elasticMM} For any $\alpha,\alpha' \in (0,1)$ and any $f \in W^{1,1}_1(m^{-1})$ and any $g \in L^1_1(m^{-1})$, there holds
$$\|\Q_\alpha^+(f,g)-\Q_{\alpha'}^+(f,g)\|_{\mathcal{X}} \leq \mathsf{p}(\alpha-\alpha')\|f\|_{W^{1,1}_1(m^{-1})}\,\|g\|_{\mathcal{Y}}$$
and
$$\|\Q_\alpha^+(g,f)- \Q_{\alpha'}^+(g,f)\|_{\mathcal{X}} \leq  \mathsf{p}(\alpha-\alpha')\|f\|_{W^{1,1}_1(m^{-1})}\,\|g\|_{\mathcal{Y}}$$
where $ \mathsf{p}(r)$ is an explicit polynomial function with $\lim_{r \to 0^+} \mathsf{p}(r)=0.$
\end{propo}
With this proposition, one sees that \eqref{XY} holds true for $\mathcal{X}=L^1(m ^{-1}).$ Moreover, arguing as in \cite[Proposition 11]{AlCaGa}, one proves easily that \eqref{XY1} holds true:
\begin{propo}\label{propoAlo} There exists $C >0$ such that, for any $\alpha  \in (0,1)$
$$\|\Q_\alpha^+(h,g)\|_{\mathcal{X}}+\|\Q_\alpha^+(g,h)\|_{\mathcal{X}} \leq C\|h\|_{\mathcal{Y}}\,\|g\|_{\mathcal{Y}} \qquad \qquad \forall h,g \in \mathcal{Y}.$$
\end{propo}
\begin{proof} The proof follows from the very simple observation that
\begin{equation}\label{Qm^1}\|\Q_\alpha^+(h,g)\|_{\mathcal{X}} \leq \|\Q_\alpha^+(m^{-1}h,m^{-1}g)\|_{L^1(\R^3)}\end{equation}
together with the well-known boundedness of the bilinear operator $\Q^+_\alpha \:: \: L^1_1(\R^3) \times L^1_1(\R^3) \to L^1(\R^3)$ (see, e.g. \cite[Theorem 1]{AlCaGa}). To prove \eqref{Qm^1}, one first notices that, for any $h,g \in \mathcal{X}$, one has
$$\|\Q_\alpha^+(h,g)\|_{\mathcal{X}}=\|\Q_\alpha^+(h,g)m^{-1}\|_{L^1}=\sup_{\|\psi\|_{L^\infty(\R^3)}=1}\int_{\R^3}\Q^+_\alpha(h,g)(v)\left(m^{-1}\psi\right)(v)\d v.$$
To estimate this last integral, one can assume without loss of generality that $h,g,\psi$ are non-negative. Then, using the weak formulation of $\Q^+_\alpha$:
$$\int_{\R^3}\Q^+_\alpha(h,g)(v)\left(m^{-1}\psi\right)(v)\d v=\int_{\R^3 \times \R^3} h(v)g(w)|v-w|\left(m^{-1}\psi\right)(v^\star)\d v\d w$$
where the post-collision velocity $v^\star$ is defined by \eqref{co:transf1}.
Now, because of the dissipation of kinetic energy, since $s \in (0,1]$, one has
$$|v^\star |^s \leq \left(|v^\star|^2+|w^\star|^2\right)^{s/2} \leq \left(|v|^2+|w|^2\right)^{s/2} \leq |v|^s+|w|^s,$$
i.e. $m^{-1}(v^\star) \leq m^{-1}(v)m^{-1}(w).$ Therefore,
$$\int_{\R^3}\Q^+_\alpha(h,g)(v)\left(m^{-1}\psi\right)(v)\d v \leq \int_{\R^3 \times \R^3} \left(m^{-1}h\right)(v)\left(m^{-1}g\right)(w)|v-w| \psi(v^\star)\d v\d w.$$
One recognizes that this last integral is equal to $\ds\int_{\R^3}\Q^+_\alpha(m^{-1}h,m^{-1}g)(v) \psi (v)\d v$ and this proves \eqref{Qm^1}.
\end{proof}
\subsection{Spectral properties of $\L$ and $\mathscr{L}_1$ in $\mathcal{X}$} The spectral properties of both the linear Boltzmann operator $\L$ and the linearized operator $\mathscr{L}$ in $\H=L^2(\M^{-1})$ are recalled in the appendix. In particular, it is well known that, for both these operators, $0$ is a simple eigenvalue associated to the eigenfunction $\M$ and both operators admit a positive spectral gap in $\H$. We shall show that the same is true in the larger space $\mathcal{X}$. To do so, we adopt the general strategy explained in the recent paper \cite{gualdani}. First, we notice that $\H$ is a dense subspace of $\mathcal{X}$. Moreover, if we denote, as in the appendix, $\mathscr{L}_2$ as the linearized Boltzmann operator $\Q_1(\cdot,\M)+\Q_1(\M,\cdot)+\L$ in $\H$, one has
$$ \mathscr{L}_{2} ={\mathscr{L}_1}_{|\H }$$
with $\mathscr{L}_1\::\:\D(\mathscr{L}_1) \subset \mathcal{X} \to \mathcal{X}$ by
$$\mathscr{L}_1 (h)=\Q_1(\M,h)+\Q_1(h,\M) +\L h , \qquad \forall h \in \D(\mathscr{L}_1)=\mathcal{Y}.$$
In the same way, with the notations of the appendix, $\L_{|\mathcal{H}}={\mathbf{L}}$. For the linear Boltzmann operator $\L$, we have the following
\begin{theo}\label{spectL} The spectrum of $\L$ in $\mathcal{X}$
  coincides with that of $\mathbf{L}$ in $\H$.  As a consequence,
  $\mathscr{N}(\L)=\mathrm{span}(\M)$ and $\L$ admits a positive
  spectral gap $\nu >0$. In particular, if
$$\widehat{\mathcal{X}}=\{f \in \mathcal{X}\,;\,\IR f \d v=0\}, \qquad \widehat{\mathcal{Y}}=\{f \in \mathcal{Y}\,;\,\IR f \d v =0\}$$ then $\mathscr{N}(\L) \cap \widehat{\mathcal{X}}=\{0\}$ and $\L$ is invertible from $\widehat{\mathcal{Y}}$ to $\widehat{\mathcal{X}}$.
\end{theo}
\begin{proof} As already mentioned, we adopt the general strategy explained in the recent paper~\cite{gualdani}.  Precisely, one proves that $\L$ splits as
$$\L=\mathcal{A}+\mathcal{B}$$
where \begin{enumerate}[(i)]\item $\mathcal{A}\::\:\mathcal{X} \to \H$ is bounded;
\item the operator $\mathcal{B}\::\:\D(\mathcal{B}) \to \mathcal{X}$ (with $\D(\mathcal{B})=\mathcal{Y}$) is $a$-dissipative for some positive $a >0$, i.e.
\begin{equation}\label{Bdiss}
\IR \mathrm{sign}f(v) \mathcal{B}f(v)m^{-1}(v)\d v \leq -a \|f\|_{\mathcal{X}} \qquad \forall f \in \mathcal{Y}.\end{equation}
\end{enumerate}
To do so, we use the estimates on $\L$ derived in the appendix. For any $R >0$, set
$$\mathcal{A}f(v)=\L^+(\chi_{\{|\cdot| \leq R\}} f)(v)=\int_{|w|\leq R} k(v,w)f(w)\d w.$$
Using Minkowski's integral inequality (with measures $\M^{-1}(v)\d v$ and $|f(w)|\d w$) one gets easily
$$\|\mathcal{A}f\|_{\H} \leq \int_{|w|\leq R} |f(w)|\left(\IR k^2(v,w)\M^{-1}(v)\d v\right)^{1/2}\d w.$$
Now, still with the notations of the appendix, one has
$$\IR k^2(v,w)\M^{-1}(v)\d v = \M^{-1}(w)\IR G^2(v,w)\d v$$
and, using Lemma A.\ref{lemmG}, with $p=2$ and $q=0$, there is some positive constant $c^2 >0$ such that
$$\IR k^2(v,w)\M^{-1}(v)\d v \leq c^2 \M^{-1}(w)(1+|w|)^{-1} \leq c^2\M^{-1}(w) \qquad \forall w \in \R^3.$$
Thus,
$$\|\mathcal{A}f\|_{\H} \leq c  \int_{|w|\leq R} |f(w)|\M^{-1/2}(w)\d w$$
and, since the domain of integration is bounded, there is some positive constant $c_R >0$ such that $\|\mathcal{A}f\|_{\H} \leq c_R \|f\|_{\mathcal{X}}$ which proves point (i). Now, we prove that $R >0$ can be chosen in such a way that
$$\mathcal{B}f(v)=\L f(v) -\mathcal{A}f(v)=\L^+(\chi_{\{|\cdot| > R\}}f)(v)-\sigma(v)f(v)$$
satisfies the above point (ii). For any $f \in \mathcal{Y}$, set $I(f)=\IR \mathrm{sign}f(v) \mathcal{B}f(v)m^{-1}(v)\d v$.
One has
\begin{equation*}\begin{split}
I(f)&=\IR \mathrm{sign}f(v) \L^+(\chi_{\{|\cdot| > R\}}f)f(v)m^{-1}(v)\d v-\IR \sigma(v)|f(v)|m^{-1}(v)\d v\\
&=\IR \mathrm{sign}f(v)m^{-1}(v)\d v \int_{\{|w|> R\}} k(v,w)\d w -\IR \sigma(v)|f(v)|m^{-1}(v)\d v  \\
&\leq \int_{\{|w|> R\}} |f(w)|H(w)\d w -\sigma_0\IR (1+|v|)|f(v)|m^{-1}(v)\d v
\end{split}
\end{equation*}
where we used the fact that $\sigma(v) \geq \sigma_0 (1+|v|)$ for some positive constant $\sigma_0 >0$ and set, as in the appendix,
$$H(w)=\IR k(v,w)m^{-1}(v)\d v, \qquad \forall w \in \R^3.$$
Then, using Proposition A.\ref{estiateH}, there is some positive constant $K >0$ such that
$$I(f) \leq K \int_{\{|w|> R\}} |f(w)| \left(1+|w|^{1-s}\right)m^{-1}(w)\d w -\sigma_0\IR (1+|v|)|f(v)|m^{-1}(v)\d v.$$
In other words,
\begin{multline*} I(f)
  \leq  - \sigma_0 \int_{\{|v| \leq R\}} |f(v)| m^{-1}(v) \,\d v
  \\
 + \int_{\{|v|> R\}} |f(v)|
  \left(
    K (1 + |v|^{1-s})
    -
    \sigma_0 (1+|v|)
  \right)
  \, m^{-1}(v)  \,\d v.\end{multline*}
We choose now $R > 0$ such that $K (1 + |v|^{1-s}) - \sigma_0
(1+|v|) \leq - \sigma_0$ for all $|v| > R$ (which can be done since $s
> 0$), so that
\begin{equation}
  \label{eq:If-3}
  I(f)
  \leq
  - \sigma_0 \IR |f(v)| m(v)^{-1} \,dv
  =
  - \sigma_0 \|f\|_{\mathcal{X}},
\end{equation}
i.e. $\mathcal{B}$ satisfies \eqref{Bdiss} with $a=\sigma_0 $. We conclude now with \cite{gualdani}.
 \end{proof}

An analogous result holds true for the linearized Boltzmann operator $\mathscr{L}$. Precisely,
\begin{theo}\label{spect}
  The spectrum of $\mathscr{L}_1$ in $\mathcal{X}$ coincides with that
  of $\mathscr{L}_2$ in $\H$ with
  $\mathscr{N}(\mathscr{L}_1)=\mathrm{span}(\M)$ and $\mathscr{L}_1$
  admits a positive spectral gap $\nu >0$. In particular,
  $\mathscr{N}(\mathscr{L}_1) \cap \widehat{\mathcal{X}}=\{0\}$ and
  $\mathscr{L}_1$ is invertible from $\widehat{\mathcal{Y}}$ to
  $\widehat{\mathcal{X}}$.
\end{theo}
\begin{proof} The proof follows exactly the same lines as the above Theorem \ref{spectL}. Precisely, one proves that $\mathscr{L}_1$ splits as
$\mathscr{L}_1=\mathscr{A}+\mathscr{B}$
where $\mathscr{A}$ and $\mathscr{B}$ are such that $\mathscr{A}\::\:\mathcal{X} \to \H$ is bounded while the operator $\mathscr{B}\::\:\D(\mathscr{B}) \to \mathcal{X}$ (with $\D(\mathscr{B})=\mathcal{Y}$) satisfies
\begin{equation}\label{Bdiss1}
\IR \mathrm{sign}f(v) \mathscr{B}f(v)m^{-1}(v)\d v \leq -a \|f\|_{\mathcal{X}} \qquad \forall f \in \mathcal{Y}\end{equation}
for some $a >0$.
One recalls that $\mathscr{L}_1f=\mathscr{L}f +\L f$ where $\mathscr{L}f=\Q_1(f,\M)+\Q_1(\M,f)$. Now, it is well known that
$$\mathscr{L}f = \mathscr{L}^+ f-\sigma_1(v) f(v) - \M(v)\left(\IR |v-w|f(w)\d w\right)$$
where  $\sigma_1(v)=\IR |v-w|\M(w)\d w \geq \sigma_1 (1+|v|)$ for some positive constant $\sigma_1 >0$ and $\mathscr{L}^+ h=\Q_1^+(h,\M)+\Q^+_1(\M,h)=2\Q^+_1(h,\M).$
Then, it is easy to recognize (see e.g. \cite{glassey}) that
$$\mathscr{L}^+ h(v)=\IR K_1(v,w)h(w)\d w$$
with
$$K_1(v,w)=C_1|v-w|^{-1} \exp
\left\{-\beta_1\left(|v-w|+\dfrac{|v|^2-{|w|}^2}{|v-w|}\right)^2\right\},$$
where $C_1 >0$ and $\beta_1=\frac{1}{8\Theta^\#}.$ In other words,  $\mathscr{L}^+$ has exactly the same form of $\L^+$ (with $\beta_0$ replaced by $\beta_1$). In particular, Proposition A. \ref{estiateH} still holds if $k(v,w)$ is replaced by $K_1(v,w)$. One defines then
\begin{equation*}\begin{split}
\mathscr{A}_1f(v)&=\mathscr{L}^+(\chi_{\{|\cdot| \leq R\}} f)(v)+ \L^+(\chi_{\{|\cdot| \leq R\}} f)(v)\\
&=\int_{|w|\leq R} \left(k(v,w)+K_1(v,w)\right)f(w)\d w.\end{split}\end{equation*}
and $\mathscr{A}_2 f(v)=-\M(v)\left(\IR |v-w|f(w)\d w\right).$ It is clear that $\mathscr{A}_2$ is bounded from $\mathcal{X}$ to $\H$ with the very rough estimate $\|\mathscr{A}_2 f\|_{\H}^2 \leq \|f\|_\mathcal{X}^2 \IR (1+|v|)^2 \M(v)\d v.$ Moreover, with the same estimates as above (using the fact that the expression of $K_1(v,w)$ is very similar to that of $k(v,w)$), one proves that, for any $R >0$, $\mathscr{A}_1\::\:\mathcal{X} \to \H$ is bounded.   We define then $\mathscr{A}=\mathscr{A}_1+\mathscr{A}_2$ so that $\mathscr{A} \::\:\mathcal{X} \to \H$ is a bounded operator. Now, set
$$\mathscr{B}f(v)=\mathscr{L}_1f(v)-\mathscr{A}f(v)=\mathscr{L}^+(\chi_{\{|\cdot| > R\}} f)(v)+ \L^+(\chi_{\{|\cdot| > R\}} f)(v)-\nu(v)f(v)$$
where $\nu(v)=\sigma(v)+\sigma_1(v) \geq \nu_* (1+|v|)$ with
$\nu_*=\sigma_0+\sigma_1.$ The estimates in the proof of the above
Theorem \ref{spectL} show then that there exists $R >0$ large enough such
that $\mathscr{B}$ satisfies~\eqref{Bdiss1} with $a=\nu_*.$ We
conclude as in \cite{gualdani}.
\end{proof}

\subsection{The elastic limit $\alpha \to 1$}\label{elasticlimit}

 The main result of this section provides a
following quantitative estimate on the distance between $\fe$ and the
Maxwellian $\M$. We recall the definition of $\mathfrak{S}_\alpha(\varrho)$ in \eqref{Salpha} and, for simplicity, $\mathfrak{S}_\alpha(1)$ shall be denoted
$\mathfrak{S}_\alpha.$ One has the following:
\begin{theo}
  \label{limit1}
  There exists an explicit function $\eta_1(\alpha)$ such that
  $\lim_{\alpha \to 1}\eta_1(\alpha)=0$ and such that for any
  $\alpha_0 \in (0,1]$
  \begin{equation*}
    \sup_{\fe \in \mathfrak{S}_\alpha}
    \left\|\fe-\M\right\|_{\mathcal{Y}}
    \leq
    \eta_1(\alpha)
    \qquad \forall \alpha \in (\alpha_0,1].
  \end{equation*}
\end{theo}

\begin{nb}
  The fact that the above conclusion does not necessarily hold for
  inelasticity parameter $\alpha \simeq 0$ is related to the
  estimate \eqref{estimDHa} hereafter. Notice that this is no major
  restriction since the above result has to be interpreted as a result
  of uniform convergence to $\M$ whenever the inelasticity parameter
  $\alpha$ goes to $1.$
\end{nb}

Let us now come to the proof of Theorem \ref{limit1} which follows the paths of the corresponding result in \cite{MiMo2}.
Let  $\M_\alpha$ denote the Maxwellian with the same mass, momentum and temperature as~$\fe$:
$$\M_\alpha(v)=\left(\dfrac{1}{2\pi\Theta_\alpha}\right)^{3/2} \exp\left(-\dfrac{|v-\ua|^2}{2\Theta_\alpha}\right)$$
where
\begin{equation}\label{uaTa}\ua=\IR v \fe(v)\d v \in \R^3 \quad \text{ and } \quad \Theta_\alpha=\dfrac{1}{3}\IR |v-\ua|^2 \fe(v)\d v > 0.\end{equation}

One can prove the following result
\begin{propo}\label{L1q} Let $\alpha_0 \in (0,1]$ be fixed. For any $q >0$ and any $\delta >0$, there is $C_\delta(q) >0$ such that the estimate
\begin{equation}\label{fel1q}\|\fe-\M_\alpha\|_{L^1_q}^{2+\delta} \leq C_\delta(q) (1-\alpha) \qquad \forall \fe \in \mathfrak{S}_\alpha\,,\forall  \alpha \in (\alpha_0,1].\end{equation}
\end{propo}
\begin{proof} Let $\alpha \in (\alpha_0,1]$ be fixed and let $\fe \in \mathfrak{S}_\alpha.$ The stationary solution $\fe$ satisfies
$$\IR \Q_\alpha(\fe,\fe)\log\bigg(\frac{\fe}{\M}\bigg)\d v =-\IR \L(\fe)\log\bigg(\frac{\fe}{\M}\bigg)\d v$$
which, from \cite[Theorem 2.1]{LoTo}, yields
$$\IR \Q_\alpha(\fe,\fe)\log\bigg(\frac{\fe}{\M}\bigg)\d v \geq 0, \qquad \forall \alpha \in (0,1].$$
Now,

$$\IR \Q_\alpha(\fe,\fe)\log\bigg(\frac{\fe}{\M}\bigg)\d v = \IR \Q_\alpha(\fe,\fe)\log \fe \d v+ \frac{1}{2\Theta^\sharp} \IR
\Q_\alpha(\fe,\fe)|v|^2 \d v$$
and, using \eqref{dissDHal}
$$\IR \Q_\alpha(\fe,\fe)\log \fe \d v=-\D_{H,\alpha}(\fe) + \frac{1-\alpha^2}{\alpha^2}\IRR \fe(v)\fe(w)|v-w|\d v \d w$$
while
$$\IR
\Q_\alpha(\fe,\fe)|v|^2 \d v=-\frac{1-\alpha^2}{8}\IRR
\fe(v)\fe(w)|v-w|^3\d v \d w.$$ Consequently, one has
 \begin{multline}
   \label{estimDHa}
  \D_{H,\alpha}(\fe) \leq \frac{1-\alpha^2}{\alpha^2}\IRR \fe(v)\fe(w)|v-w|\d v \d w \\ -\frac{1-\alpha^2}{16\Theta^\sharp}\IRR \fe(v)\fe(w)|v-w|^3\d v \d w \\
  \leq \frac{1-\alpha^2}{\alpha^2}\IRR \fe(v)\fe(w)|v-w|\d v \d w
\end{multline}
From the estimate of the moments of $\fe$, this last integral can be
estimated from above by some positive constant $K >0$ independent of
$\alpha \in (0,1]$. In particular, for any \textit{fixed} $\alpha_0
\in (0,1]$, there is $C_0 >0$ such that
$$\D_{H,\alpha}(\fe) \leq C_0 (1-\alpha) \qquad \forall \alpha \in (\alpha_0,1].$$
The above estimate, together with Proposition \ref{estimentropy}, implies the existence of some $C_1 >0$ such that
\begin{equation}\label{dh1}\D_{H,1}(\fe) \leq  C_1  (1-\alpha) \qquad \forall \alpha \in (\alpha_0,1].\end{equation}
Recall that $\D_{H,1}$ is the entropy dissipation functional associated to classical (elastic) interactions and has been studied intensively in \cite{villcerc}. In particular, using the estimates of Corollary \ref{sobolev} and Theorem \ref{inf2t}, one deduces from the \textit{op. cit.} (see also \cite[Theorem 3.5]{MiMo3}) that, for any $\delta >0$, there is $\widetilde{C}_\delta >0$ such
\begin{equation}\label{villa}\|\fe-\M_\alpha\|_{L^1}^{2} \leq 2  \IR \fe(v)\log \frac{\fe(v)}{\M_\alpha (v)} \d v \leq \widetilde{C}_\delta \D_{H,1}(\fe)^{\frac{2}{2+\delta}}.\end{equation}
Then, from \eqref{dh1}, we get that, for any $\delta >0$, there exists $C_\delta >0$ such that
$$\|\fe-\M_\alpha\|_{L^1} \leq C_\delta (1-\alpha)^{\frac{1}{1+\delta}} \qquad \forall \alpha \in (\alpha_0,1).$$
Now, using Theorems \ref{gaus} and \ref{theoenergy}, by a simple interpolation argument, we get the conclusion. \end{proof}

An easy consequence of the above Theorem is the following where we recall that the space $\mathcal{Y}$ has been defined in the previous section.
\begin{cor} \label{FaMaX} For any $\delta >0$,   there exists an explicit constant  $C_\delta  >0$  such that, for any
  $\alpha_0 \in (0,1]$
  \begin{equation*}
    \left\|\fe-\M_\alpha\right\|_{\mathcal{Y}}
    \leq
    C(1-\alpha)^{\frac{1}{4+2\delta}}
    \qquad \forall \fe \in \mathfrak{S}_\alpha\,,\quad \forall \alpha \in (\alpha_0,1].
  \end{equation*}
 \end{cor}
\begin{proof} The proof relies on a simple interpolation argument from Proposition \ref{L1q} (with $q=1$) and Theorem \ref{gaus}. Recall that $\mathcal{Y}=L^1_1(m^{-1})$ where $m(v)=\exp(-a|v|^s)$ for some fixed $a >0$, $s \in (0,2)$. For any $\alpha\in (\alpha_0,1]$ and any $\fe \in \mathfrak{S}_\alpha$, one has
\begin{multline*} \left\|\fe-\M_\alpha\right\|_{\mathcal{Y}} \leq \left(\IR \left|\fe(v)-\M_\alpha(v)\right| \langle v \rangle \d v \right)^{1/2}\\ \left(\IR \left|\fe(v)-\M_\alpha(v)\right| \langle v \rangle  \exp(2a|v|^s)\d v\right)^{1/2}.
\end{multline*}
Moreover, according to Theorem \ref{gaus} (and since the energy $E_\alpha$ of $\fe$ can be bounded from below and above independently of $\alpha \in (0,1)$), there exist  $A >0$ and $M >0$ such that
$$\IR \left|\fe(v)-\M_\alpha(v)\right| \exp(A|v|^2)\d v \leq M \qquad \forall \alpha \in (\alpha_0,1].$$
Since there exists $c=c(a,s)$ such that $\langle v \rangle  \exp(2a|v|^s) \leq c(a,q,s) \exp(A|v|^2)$ for any $v \in \R^3$, one gets the conclusion with $C=\sqrt{M\,c(a,s)C_\delta(1)}$ where $C_\delta(1)$ is the constant appearing in Prop. \ref{L1q}.
\end{proof}

With the above Corollary, one gets the following:
\begin{lemme}  There exist explicit constants $C >0$ and $p >0$ such that, for any   $\alpha_0 \in (0,1]$,
$$\|\L(F_\alpha)\|_{\mathcal{X}}=\|\Q_\alpha(\fe,\fe)\|_{\mathcal{X}} \leq C(1-\alpha)^p \qquad \forall \fe \in \mathfrak{S}_\alpha\,,\:\:\forall \alpha \in (\alpha_0,1].$$
\end{lemme}
\begin{proof} Let $\alpha_0 \in (0,1]$ be fixed and let $\alpha \in (\alpha_0,1]$. For any $\fe \in \mathfrak{S}_\alpha$, one has
\begin{equation*}\begin{split}
-\L(\fe)&=\Q_\alpha(\fe,\fe)=\Q_\alpha(\fe,\fe-\M_\alpha)+\Q_\alpha(\fe,\M_\alpha)\\
&=\Q_\alpha(\fe,\fe-\M_\alpha)+\Q_\alpha(\fe-\M_\alpha,\M_\alpha)+\Q_\alpha(\M_\alpha,\M_\alpha).
\end{split}\end{equation*}
Thus,
$$\|\L(\fe)\|_{\mathcal{X}} \leq \|\Q_\alpha(\fe,\fe-\M_\alpha)\|_{\mathcal{X}}+\|\Q_\alpha(\fe-\M_\alpha,\M_\alpha)\|_{\mathcal{X}}+\|\Q_\alpha(\M_\alpha,\M_\alpha)\|_{\mathcal{X}}.$$
Then, from Prop. \ref{propoAlo}, there exists $C >0$ such that
$$\|\L(\fe)\|_{\mathcal{X}} \leq C\,\|\fe-\M_\alpha\|_{\mathcal{Y}}\left(\|\fe\|_{\mathcal{Y}}+\|\M_\alpha\|_{\mathcal{Y}}\right)+\|\Q_\alpha(\M_\alpha,\M_\alpha)\|_{\mathcal{X}}.$$
Moreover, since $\M_\alpha$ is the Maxwellian with same first moments as $\fe$, it is easy to see that $\|\M_\alpha\|_{\mathcal{Y}}$ depends only on the energy $E_\alpha=\IR\fe(v)|v|^2\d v$. Thus, on the basis of the a posteriori estimates derived in Section 2, namely Theorem \ref{theoenergy}, one gets easily that
$$\sup_{\alpha \in (0,1]}\left(\|\fe\|_{\mathcal{Y}}+\|\M_\alpha\|_{\mathcal{Y}}\right) < \infty.$$
Therefore, there exists a positive constant $C_2 >0$ such that
$$\|\L(\fe)\|_{\mathcal{X}} \leq C_2\,\|\fe-\M_\alpha\|_{\mathcal{Y}}+\|\Q_\alpha(\M_\alpha,\M_\alpha)\|_{\mathcal{X}}.$$
Now, to estimate $\|\Q_\alpha(\M_\alpha,\M_\alpha)\|_{\mathcal{X}}$, one only notices that, since $\M_\alpha$ is a Maxwellian, one has $\Q_1(\M_\alpha,\M_\alpha)=0$, i.e.
$$\|\Q_\alpha(\M_\alpha,\M_\alpha)\|_{\mathcal{X}}=\|\Q_\alpha(\M_\alpha,\M_\alpha)-\Q_1(\M_\alpha,\M_\alpha)\|_{\mathcal{X}}.$$
Therefore, one can apply Proposition \ref{elasticMM} to get the existence of some polynomial mapping $r \mapsto \mathsf{p}(r)$ such that
$$\|\Q_\alpha(\M_\alpha,\M_\alpha)\|_{\mathcal{X}} \leq \mathsf{p}(1-\alpha)\|\M_\alpha\|_{W^{1,1}_1(m^{-1})}\,\|\M_\alpha\|_{\mathcal{X}}$$
where $\lim_{r \to 0}\mathsf{p}(r)=0.$ Again, since the various norms of $\M_\alpha$ only depend on the energy~$E_\alpha,$ we deduce from Theorem \ref{theoenergy} that there exists a positive constant $C_3$ such that
$$\|\Q_\alpha(\M_\alpha,\M_\alpha)\|_{\mathcal{X}} \leq C_3\mathsf{p}(1-\alpha) \qquad \forall \alpha \in (\alpha_0,1].$$
Consequently, there exist two positive constants $C_2,C_3 >0$ and some polynomial function $r \mapsto \mathsf{p}(r)$ with $\lim_{r \to 0}\mathsf{p}(r)=0$ such that
$$\|\L(\fe)\|_{\mathcal{X}}\leq C_2\,\|\fe-\M_\alpha\|_{\mathcal{Y}}+C_3\mathsf{p}(1-\alpha) \qquad \forall \alpha \in (\alpha_0,1].$$
We get the desired estimate using Corollary \ref{FaMaX}.
\end{proof}
The above Lemma allows us to conclude the proof of Theorem \ref{limit1}:
\begin{proof}[Proof of Theorem \ref{limit1}] For any $\alpha \in (0,1)$, set  $g_\alpha=-\Q_\alpha(\fe,\fe)$, we get
$$\L(\fe)=\L(\fe-\M)=g_\alpha$$
with $\fe-\M \in \widehat{\mathcal{Y}}.$ Since $\L$ is invertible from $\widehat{\mathcal{Y}}$ to $\widehat{\mathcal{X}}$ (with bounded inverse) according to Theorem \ref{spectL}, there is some positive constant $c >0$ such that
$$\|\fe-\M\|_{\mathcal{Y}}=\|\L^{-1}(g_\alpha)\|_{\mathcal{Y}}\leq c\|g_\alpha\|_{\mathcal{X}}.$$
According to the above Lemma $\lim_{\alpha \to 1}\|g_\alpha\|_{\mathcal{X}}=0$ which yields the result.\end{proof}
\begin{nb} Notice that Theorem \ref{limit1} combined with Corollary \ref{FaMaX} shows that
$$\lim_{\alpha \to 1}\|\M_\alpha -\M\|_{\mathcal{Y}}=0$$
with some explicit rate, where $\M_\alpha$ is the Maxwellian with same
mass, momentum and temperature as~$\fe$. This implies in particular
that
$$\lim_{\alpha \to 1} \ua=0 \qquad \text{ while } \qquad \lim_{\alpha \to 1} \Theta_\alpha=\Theta^\#$$
where $\ua$ and $\Theta_\alpha$ are defined in \eqref{uaTa}. Since
$\L$ does not conserve momentum, it is not clear how to prove
convergence of the first moments of $\fe$ towards those of $\M$ in a
direct way. Notice that, for the forcing terms considered in previous
related works \cite{MiMo2,MiMo3}, the convergence of $\fe$ to $\M$
was, on the contrary, proved thanks to the convergence of the momentum
and temperature.
\end{nb}
\subsection{Uniqueness} With this in hands, as explained at the beginning of the section, one can state the following:
\begin{theo}\label{uniqueness} There exists $\alpha_0 \in (0,1]$ such that, for any $\varrho > 0$, the set
$$\mathfrak{S}_\alpha(\varrho)=\left\{F_\alpha \in L^1_2, \:F_\alpha \geq 0,\; \fe  \text{ solution to \eqref{equi} with } \IR \fe \d v =\varrho  \right\}$$
reduces to a singleton. In particular, for any $\alpha \in (\alpha_0,1]$, such a steady state $F_\alpha$ is radially symmetric and belongs to $\mathcal{C}^\infty(\R^3)$.
\end{theo}
\begin{proof} Our strategy to prove the uniqueness result has been explained in Section \ref{sec:uniqueness} and we have already shown that the estimates \eqref{XY}, \eqref{XY1} and \eqref{XY2} hold true with our choice of $\mathcal{X}$ and $\mathcal{Y}$ while condition  \eqref{deltac0} holds thanks to Theorem \ref{limit1}.
\end{proof}
\section*{Appendix: Properties of the linear and linearized Boltzmann operators}\setcounter{equation}{0}
\renewcommand{\theequation}{A.\arabic{equation}}

\subsection*{A.1. Spectral analysis of the linearized operator in $L^2(\M^{-1})$} We consider here the spectral analysis of the linearization of $\B(f,f)+\L(f)$ around the Maxwellian state $\M$. Precisely, let $\H$ denote the Hilbert space $L^2(\M^{-1})$ endowed with the inner space
$$\langle f,g\rangle_\H :=\IR f(v)g(v)\M^{-1}(v)\d v, \qquad \forall f,g \in \H$$
and let $\LL$ denote the following unbounded operator in $\H$:
$$\LL(h)=\mathbf{L}(h) +\L(h), \qquad \forall h \in \mathscr{D}(\LL)$$
where $\mathbf{L}(h)=\Q_1(h,\M)+\Q_1(\M,h)$ is the linearized operator of the classical Boltzmann operator $\Q_1(\cdot,\cdot)$. The domain of $\LL$ in $\H$ is
$$\mathscr{D}(\LL)=L^2_1(\M^{-1})=\left\{f=f(v)\;;\;\IR |f(v)|^2  \M^{-1}(v)\left(1+|v|^2\right)^{1/2}\d v < \infty\right\}.$$
The spectral analysis of the linearized operator $\mathbf{L}$ in $\H$ is a well-known feature of the classical theory Boltzmann operator (see \cite[Chapter 7]{Ce94}, \cite[Chapter 3]{glassey}). In particular, $\mathbf{L}$ is a nonnegative self-adjoint operator in $\H$ with
\begin{equation}\label{a1}\langle h, \mathbf{L}h\rangle_\H \leq 0 \qquad \forall h \in \mathscr{D}(\LL)\end{equation}
and $N(\mathbf{L})=\mathrm{span}\left\{\M,v_1\,\M,v_2\,\M,v_3\,\M,|v|^2\,\M\right\}.$ Moreover, the spectral analysis of $\L$ in $\H$ has been performed in \cite{ArLo} and made precise in \cite{LoMo}. Here again, $\L$ is a nonnegative self-adjoint operator in $\H$ and  there exists $\mu  >0$ such that
\begin{equation}\label{a2}-\langle f,\L(f)\rangle_\H \geq \mu \|f-\varrho_f \mathcal{M}\|_{L^2(\mathcal{M}^{-1})}^2 \qquad \forall f \in \mathscr{D}(\LL)\end{equation}
with $\varrho_f =\IR f(v)\d v$ and where some quantitative estimates of the spectral gap $\mu$ are given in \cite{LoMo}.
In particular, $N(\L)=\mathrm{span}\left\{\M \right\}$. One deduces directly from \eqref{a1} and~\eqref{a2} that
$$-\langle f, \LL(f)\rangle_\H \geq  \mu \|f-\varrho_f \mathcal{M}\|_{L^2(\mathcal{M}^{-1})}^2 \qquad \forall f \in \mathscr{D}(\LL).$$
In particular,
$$N(\LL)=\mathrm{span}\left\{\M \right\}.$$
Moreover, it is not difficult to resume the arguments of both \cite{glassey} and \cite{ArLo} to prove that there exists some nonnegative measurable function $\nu(v)$ such that
$$\LL(f)=\LL^c(f)- \nu(v) f(v)$$ where $\LL^c$ is an integral operator, relatively compact with respect to the multiplication operator $f \mapsto \nu f$. Therefore, the spectrum $\mathfrak{S}(\LL)$ of $\LL$ is made of continuous (essential) spectrum
$\{\lambda \in \mathbb{R}\,;\, \lambda \leq -\nu_0\}$ where
$\nu_0=\inf_{\v \in \R^3}\nu(\v)
>0$ and a  decreasing sequence of real eigenvalues with finite algebraic
multiplicities whose unique possible cluster point is  $-\nu_0.$ Moreover, the spectral gap $\mu_2$ of $\LL$
$$\mu_2:=\min\bigg\{\lambda  :  -\lambda \in (-\nu_0,0), - \lambda \in
\mathfrak{S}(\LL)\setminus \{0\}\bigg\}$$
satisfies the quantitative
estimate $\mu_2 \geq \mu$ where $\mu$ is the spectral gap of $\L$
given in \eqref{a2}.

\subsection*{A.2. Estimates on the linear operator $\L$}
\setcounter{equation}{0}

We now establish several important estimates on the linear Boltzmann operator $\L$. Precisely, we recall first the spectral properties of $\L$ in~$\mathcal{H}$
where we recall that $\mathcal{H}=L^2(\M^{-1}\d v)$. To distinguish the linear Boltzmann operator in $\H$ and in $\mathcal{X}$, one shall denote by $\mathbf{L}$ the linear Boltzmann operator in $\H$: $\mathbf{L}\::\:\D(\mathbf{L}) \subset \H \to \H$ with
$$\mathscr{D}(\mathbf{L})=L^2_1(\M^{-1})=\left\{f=f(v)\;;\;\IR |f(v)|^2  \M^{-1}(v)\left(1+|v|^2\right)^{1/2}\d v < \infty\right\}$$
and,
$$\mathbf{L} f(v)=\int_{\R^3}k(v,w)f(w)\d w-\sigma(v)f(v), \qquad \forall f \in \mathscr{D}(\mathbf{L})$$
where $k(v,w)$ is given by \eqref{k} and $\sigma(\cdot)$ is defined in \eqref{sigma} and satisfies
$$\sigma(v)=\int_{\R^3}k(v,w)\d w \geq \sigma_0(1+|v|) \qquad \forall v \in \R^3$$
 with $\sigma_0 >0$. Moreover, the spectral structure of $\mathbf{L}$ has been studied in \cite{ArLo, LoMo} and can be summarized in the following:
\begin{propoA} The spectrum $\mathfrak{S}(\mathbf{L})$ of the operator
$\mathbf{L}$  in $\H$ is made of continuous (essential) spectrum
$\{\lambda \in \mathbb{R}\,;\, \lambda \leq -\nu_0\}$ where
$\nu_0=\inf_{\v \in \R^3}\sigma(\v)
>0$ and a  decreasing sequence of real eigenvalues with finite algebraic
multiplicities whose unique possible cluster point is  $-\nu_0$. Moreover, $0$ is an eigenvalue of $\mathbf{L}$
associated to $\M$ and $\mathbf{L}$ admits a  spectral gap $\mu_0 >0$ such that
\begin{equation*}
\mu_0 :=\min\bigg\{\lambda  :  -\lambda \in (-\nu_0,0), - \lambda \in
\mathfrak{S}(\mathbf{L})\setminus \{0\}\bigg\}\geq \frac{\eta(1+e)}{4\sqrt{5}} >0
\end{equation*}
with  $\eta=\sqrt{2\Theta_0}
\, \mathrm{erf}^{-1}\left(\frac{1}{2}\right)$ where $\mathrm{erf}^{-1}$
denotes the inverse error function,
$\mathrm{erf}^{-1}(\frac{1}{2})\simeq 0.4769.$
\end{propoA}

Notice that several properties of the kernel $k(v,w)$ have been derived in \cite{ArLo} in the spirit of \cite{carleman}. Precisely, one has
$$k(v,w)\M(w)=k(w,v)\M(v) \qquad \forall v,w \in \R^3 \times \R^3$$
and, setting
$$G(v,w)=\M^{-1/2}(\v)k(\v,w){\M}^{1/2}(w),\qquad \qquad \v,w \in \R^3 \times \R^3,$$
one has $G(v,w)=G(w,v)$ and the following holds
\begin{lemmeA}\label{lemmG}
For any $0 < p < 3$ and any $q \geq 0$, there exists $C(p,q) >0$
such that
$$\int_{\R^3}|G(v,w)|^p \dfrac{\d v}{(1+|v|)^q} \leq \dfrac{C(p,q)}{(1+|w|)^{q+1}}, \qquad \forall w \in \R^3.$$
\end{lemmeA}
We shall exploit this estimate to derive the following more general one in which algebraic weights are replaced by exponential weight. Namely, one proves the following:
\begin{propoA}\label{estiateH} Set $m(v)=\exp(-a|v|^s)$, $a >0$ and $s \in (0,1]$ and
$$H(w)=\int_{\R^3} k(v,w)m^{-1}(v)\d v, \qquad \qquad w \in \R^3.$$
Then, there exists a positive constant $K=K(e,a,s) >0$ such that
$$H(w) \leq K(1 + |w|^{1-s})\,m^{-1}(w) \qquad \forall w \in \R^3.$$
\end{propoA}
\begin{proof} Recall that $k(v,w)$ is given by \eqref{k}. Taking into account that
$\frac{|v|^2-{|w|}^2}{|v-w|} - |v-w|
  =
  2 \frac{v-w}{|v-w|}\cdot w$ we may rewrite \eqref{k} as
\begin{equation}
  \label{eq:k2}
  k(v,w)
  = C_0 |v-w|^{-1} \exp \left\{
    -\beta_0
    \left(
      (2+\mu)|v-w| + 2\frac{v-w}{|v-w|} \cdot w
    \right)^2
  \right\}.
\end{equation}
Performing the change of variables $u=v-w$ and using spherical coordinates (with $\varrho=|u|$ and $\varrho|w|y=u \cdot w$)  one gets easily
\begin{equation*}
H(w)=2\pi C_0 \int_{A}F(\varrho,y)\d\varrho\d y
\end{equation*}
with $A=[0,\infty) \times [-1,1]$ and
$$F(\varrho,y)=\varrho\exp \left\{
      -\beta_0
      \bigg(
        (2+\mu)\varrho + 2 |w|y
      \bigg)^2
      +  a       \bigg(\varrho^2+|w|^2 + 2\varrho|w|y\bigg)^{s/2}\right\}$$
Split $A$ into the two regions of integration:
$$A_1= \{(\varrho,y) \in A \, ;\,3 |w| y \geq
    - 2 \varrho \} \qquad \text{ and } \qquad A_2=A\setminus A_1.$$
Notice first that, since $y \leq 1$ and $s\in (0,1]$
$$\exp\left(a(\varrho^2+|w|^2+2\varrho|w|y)^{s/2}\right) \leq \exp\left(a(\varrho+|w|)^{s}\right) \leq \exp(a \varrho^s)\exp(a |w|^s) \qquad \forall (\varrho,y) \in A.$$
Moreover, since $(2+\mu)\varrho + 2 |w|y \geq (\mu + 2/3) \varrho$ for any $(\varrho,y) \in A_1$ we have
  \begin{equation}\begin{split}
    \label{eq:region1}
    \int_{A_1} F(\varrho,y) \d\varrho\d y
    &\leq \exp(a |w|^s) \int_0^\infty \d\varrho \! \int_{-1}^1 \varrho \exp
    \left( -\beta_0 (2/3+\mu)^2\varrho^2 + a \varrho^s \right) \,\d y
    \\
  &\leq C_1 \exp(a |w|^s) = C_1 m^{-1}(w) \end{split}\end{equation}
  since the integral is convergent.\\

Let us estimate now the integral over $A_2$ which is more intricate.  For any $(\varrho,y) \in A_2$, one notices first that
$$\varrho^2 + |w|^2 + 2 \varrho|w|y < |w|^2 -\varrho^2/3 \qquad \text{ and } \qquad \varrho \leq (3/2) |w|,$$
so that
  \begin{multline}
    \label{eq:region2}
    \int_{A_2} F(\varrho,y) \,\d y \,\d\varrho
    \\
    \leq \int_0^{(3/2)|w|} \varrho \exp \left( a \left(|w|^2 -\varrho^2/3
      \right)^{s/2} \right) \d\varrho\int_{-1}^1 \exp \left( -\beta_0 \left(
        (2+\mu)\varrho + 2 |w|y \right)^2 \right) \,\d y.
  \end{multline}
  To carry out the $y$-integral, perform the change of variables  $z =
  (2+\mu)\varrho + 2 |w|y$ to get
  \begin{equation*}
    \int_{-1}^1
    \exp \left( -\beta_0
      \left(
        (2+\mu)\varrho + 2 |w|y
      \right)^2
    \right)
    \,\d y
    \leq
    \frac{1}{2|w|}
    \int_{-\infty}^{\infty}
    \exp \left( -\beta_0 z^2
    \right)
    \,\d z
    =
    \frac{C_2}{|w|}
  \end{equation*}
  for some explicit $C_2 >0.$ Plugging this  in \eqref{eq:region2} we obtain
$$\int_{A_2} F(\varrho,y) \,\d y \,\d\varrho
    \leq \frac{C_2}{|w|}\int_0^{(3/2)|w|} \varrho \exp \left( a \left(|w|^2 -\varrho^2/3
      \right)^{s/2} \right) \d\varrho.$$
Setting now  $x=|w|^2 - \varrho^2/3$, we obtain
\begin{equation}\label{eq:region2-2}\int_{A_2}F(\varrho,y)\d\varrho\d y \leq \dfrac{3C_2}{2\,|w|}\int_{|w|^2/4}^{|w|^2}\exp(a x^{s/2})\d x \leq \dfrac{3C_2}{2\,|w|}\int_{0}^{|w|^2}\exp(a x^{s/2})\d x.\end{equation}
We observe now that, for any $r >0$, \begin{equation*}\begin{split}
    \int_0^r \exp ( a x^{s/2} ) \,\d x
    &\leq \frac{2}{as}
    r^{1-s/2} \int_0^r \frac{as}{2} x^{s/2-1}\exp ( a x^{s/2}
    ) \,\d x
    \\
    &= \frac{2}{as} r^{1-s/2} \int_0^r \frac{\d}{\d x}\exp  ( a
      x^{s/2} ) \,\d x \leq \frac{2}{as} r^{1-s/2} \exp ( a
      r^{s/2} ).\end{split}\end{equation*}
 Using this in \eqref{eq:region2-2} for $r = |w|^2$ we get
\begin{equation}
    \label{eq:region2-3}
    \int_{A_2}F(\varrho,y)\d\varrho\d y
    \leq
    \frac{3C_2}{a\,s}
    |w|^{1-s}
    \exp \left( a |w|^{s} \right).
  \end{equation}
    Putting together \eqref{eq:region1} and \eqref{eq:region2-3} we finally obtain the result.\end{proof}
\begin{nbA}
Notice  that, whenever $s=1$, the above Proposition actually asserts that $H(w) \leq Cm^{-1}(w)$ for any $w \in \R^3$. Moreover,  for any $f \in \mathcal{X}=L^1(\R^3,m^{-1}(v)\d v),$ one has
$$\|\L^+ f\|_{\mathcal{X}} \leq \int_{\R^3} |f(w)|H(w)\d w$$
where $\L^+$ is the restriction of $\mathbf{L}$ to $\mathcal{X}$. In
other words, for $s=1$, we get that $\L^+\::\:\mathcal{X} \to
\mathcal{X}$ is a bounded operator. This is reminiscent from
\cite[Theorem 12]{AlCaGa} where exponential moment estimates for
$\Q^+_e(f,g)$ (with non Maxwellian weights) are derived. Notice that,
in \cite[Theorem~12]{AlCaGa}, an assumption of strict inelasticity
(corresponding here to $e<1$) was required which is not needed in the
above Proposition.
\end{nbA}

\end{document}